\documentclass[11pt]{article}
\usepackage{amsmath,amsthm}

\def\ftoday{le \space\number\day \space\ifcase\month\or
  janvier\or f\'evrier\or mars\or avril\or mai\or juin\or
  juillet\or ao\^ut\or septembre\or octobre\or novembre\or d\'ecembre\fi
  \space\number\year}

\def\real{I\kern-0.20em R}

\def\integer{I\kern-0.20em N}
\def\relative{{\rm \rlap Z\kern 2.2pt Z}}
\def\cc{\kern-.25em{\c c}}

\def\bc{\begin{center}}
\def\ec{\end{center}}
\def\=def{\stackrel{{\rm def}}{=}}

\newcounter{indconst}
\newcounter{auxconst}

\def\bit{\begin{itemize}}
\def\eit{\end{itemize}}
\def\ben{\begin{enumerate}}
\def\een{\end{enumerate}}
\def\bde{\begin{description}}
\def\ede{\end{description}}
 
\def\beq{\begin{equation}}
\def\eeq{\end{equation}}
\def\bfi{\begin{figure}[hbt] \begin{center}}
\def\efi{\end{center} \end{figure}}

\def\bce{\begin{center}}
\def\ece{\end{center}}

\newtheorem {theos} {Theorem}[section]

\newtheorem {coros} {Corollary} [section]

\newtheorem {lemms} {Lemma}[section]
\newtheorem {prop} {Proposition}[section]

\newtheorem {defis} {Definition}[section]

\newtheorem {rems}{Remark} [section]
\newtheorem {ex} {Exemple}[section]

\input amssym.def
\input amssym
\begin{document}

%\makeindex
\title{Family of intersecting totally real manifolds of $(\Bbb C^n,0)$ and CR-singularities}
\author{Laurent Stolovitch \thanks{CNRS UMR 5580, Laboratoire Emile Picard,
Universite Paul Sabatier, 118 route de Narbonne,
31062 Toulouse cedex 4, France. Courriel : {\tt stolo@picard.ups-tlse.fr}} }
\date{\ftoday}
\maketitle
\def\abstractname{Abstract}
\begin{abstract}
The first part of this article is devoted to the study families of totally real intersecting $n$-submanifolds
of $(\Bbb C^n,0)$. We give some conditions which allow to straighten
holomorphically the family. If this is not possible to do it formally, we construct a germ of complex analytic set at the
origin which interesection with the family can be holomorphically staightened. The second
part is devoted to the study real analytic $(n+r)$-submanifolds of $(\Bbb C^n,0)$ having a
CR-singularity at the origin ($r$ is a nonnegative integer). We consider deformations of quadrics and we
define generalized Bishop invariants. Such a quadric intersects the complex
linear manifold ${z_{p+1}=\cdots=z_n=0}$ along some real linear set ${\cal L}$. 
We study what happens to this intersection when the quadric is analytically perturbed. 
On the other hand, we show, under some assumptions, that if such a submanifold is formally
equivalent to its associated quadric then it is holomorphically equivalent to
it. 

All these results rely on a result stating the existence (and
caracterization) of a germ of
complex analytic set left invariant by an abelian group of
germs of holomorphic diffeomorphisms (not tangent to the identity at the origin).
\end{abstract}
\tableofcontents
\section{Introduction}

The aim of the article is to study the geometry of some germs of real
analytic submanifolds of $(\Bbb C^n,0)$. On
the one hand, we shall study family of totally real submanifolds of $(\Bbb C^n,0)$
intersecting at the origin. On the other hand, we shall study submanifolds
having a CR-singularity at the origin. In both cases, we are primary
interested in the holomorphic classification of such objects, that is the
orbit of the action of the group of germs of holomorphic diffeomorphisms
fixing the origin. 

In this article, we shall mainly focus on the existence of complex analytic
subsets intersecting such germs of real analytic manifolds. In the first
problem, we shall also be interested in the problem of straightening
holomorphically the family. We mean that we shall give sufficient condition which will ensure that, in
a good holomorphic coordinates system, each submanifold of the family is an $n$-plane. In the
case there are formal obstructions to straighten the family, we show the
existence of a germ complex analytic variety which interesects the family along a
set that can be straightened.
The first part of this work takes its roots in and generalizes a recent work
of Sidney Webster \cite{webster-inter} from which it is very inspired. This part of the work start after having listen to Sidney Webster at the Partial Differential Equations and Several Complex
Variables conference held in Wuhan University in June 2004.

The starting point of the first problem appeared already in the work of E. Kasner \cite{kasner1} and was studied,
from the formal view point, by G.A. Pfeiffer \cite{Pfeiffer1}. They were
interested in pairs of real analytic curves in $(\Bbb C,0)$ passing throught
the origin. We shall not consider the case were some of the submanifolds
are tangent to some others. We refer the reader to the works of I. Nakai \cite{nakai-toulouse},
J.-M. Tr\'epreau \cite{trepreau-tg} and X. Gong \cite{gong-ahern-tg} in this direction.

In the second part, we shall study $(n+p-1)$-real analytic submanifolds of $\Bbb C^n$ of the form
\begin{equation*}
\begin{cases}
y_{p+1} = F_{p+1}(z',\bar z', x'')\\ \vdots\\ y_{p+q} = F_{p+q}(z',\bar z', x'')\\
z_{n} = G(z',\bar z', x'')\end{cases}
\end{equation*}
where $z'=(z_1,\ldots,z_p)\in\Bbb C^p$, $z''=(z_{p+1},\ldots, z_{n-1})\in \Bbb
C^{n-1-p}$, $p\in \Bbb N^*$. The $F$'s and $G$ vanish at the second order at the
origin. The origin is a singularity for the Cauchy-Riemann structures. 
The most studied case, up to now, is the case where $p=1$ ($n$-submanifold). Nevertheless, some work has been done for smaller dimensional submanifold by Adam Coffman (see for instance \cite{coffman-houston, coffman-cross}). 
A nondegenerate real analytic surface in $(\Bbb C^2,0)$ which is totally real except at the origin where it has a complex tangent can be regarded as a third order analytic deformation of the quadric
$$
z_2= z_1\bar z_1+ \gamma (z_1^2+ \bar{z_1}^2).
$$
This is due to Bishop \cite{bishop} and the nonnegative number $\gamma$ is called the Bishop invariant.
When $0<\gamma <1/2$ (we say elliptic), J. Moser and S. Webster showed, in their pioneering work  \cite{moser-webster} that an analytic deformation of such a quadric could be transformed into a normal form (in fact a real algebraic variety) by the mean of a germ of holomorphic diffeomorphism preserving the origin. All the geometry a such a deformation is understood by the study of the normal form.
When $1/2<\gamma$ (we say hyperbolic), it is known that such a statement doesn't hold. So, what about the geometry ? Wilhelm Klingenberg Jr. showed \cite{klingenberg} that there exists a germ of complex curve passing throught the origin cutting the submanifold along two transversal real analytic curves  which are tangent to $\{\bar z_1 =\lambda z_1\}$ and $\{\bar z_1 =\lambda^{-1} z_1\}$ respectively. 
Here, $\lambda$ is a solution of $\gamma\lambda^2+\lambda+\gamma=0$ which is assumed to satify a {\it diophantine condition}~: there exists $M,\delta>0$ such that, for all positive integer $k$, $|\lambda^k-1|>Mk^{-\delta}$. We refer to \cite{ber-bullams} for a summary in this framework,  to \cite{huang-lecture} for a nice introduction and to \cite{harris-spm} for another point of view.

We shall work with nondegenerate submanifolds (in some sense) and then
define generalized Bishop invariants $\{\gamma_i\}_{i=1,\ldots,p}$~: 
There is a good holomorphic coordinates system at the origin in which the submanifold is defined by
\begin{equation*}
(M_{Q})\begin{cases}
y_{\alpha}  =   f_{\alpha}(z',\bar z', x'')\quad\alpha=p+1,\ldots, n-1\\ 
z_{n}  =  Q(z',\bar z')+ g(z',\bar z', x'')
\end{cases}
\end{equation*}
where the $f_i$'s and $g$ are germs of real analytic functions at the origin
and of order greater than or equal to $3$ there. The quadratic polynomial $Q$ is of the form 
$$
Q(z',\bar z') = \sum d_{i,l} z'_i \bar z'_l+ \sum_{i=1}^p  \gamma_{i}( (z'_i)^2
+(\bar z'_i)^2),
$$
the norm of the sesquilinear part of $Q$ being $1$. It is regarded as a perturbation of the quadric
\begin{equation*}
(Q)\begin{cases}
y_{p+1} =  \cdots = y_{n-1} = 0\\
z_{n} = Q(z',\bar z')\end{cases}
\end{equation*}

We shall consider the case where none of these invariants vanishes (see \cite{moser-zero,huang-krantz} for results in this situation for $p=1$).
Then, we begin a study {\it \`a la Moser-Webster} of such an object althought we
shall not study  here the normal form problem. We associate a pair $(\tau_1,\tau_2)$ of germs of holomorphic involutions of $(\Bbb C^{n+p-1},0)$ and also a germ of biholomorphism $\Phi=\tau_1\circ\tau_2$. 

In the one hand, we shall show under some asumptions (in particular the diophantiness property of $D\Phi(0)$), that if $(M_Q)$ is formally equivalent to its associated quadric $(Q)$ then it is holomorphically equivalent to it. This generalizes a result by X. Gong \cite{gong-hyperbolic}(p=1). 

On the other hand, the complex linear space $\{z_{p+1}=\cdots=z_n=0\}$ intersects the quadric along the real linear set ${\cal L}:=\{Q(z',\bar z')=0\}$. We shall show that this situation survive under a small analytic perturbation. Namely, we shall show under some assuptions, that there exists a good holomorphic coordinates system at the origin in which the submanifold $(M_{Q})$ intersects the complex linear space $\{z_{p+1}=\cdots=z_n=0\}$ along a real analytic subset $V$ passing throught the origin. The latter is completely determined by the eigenvalues of $D\Phi(0)$ and more precisely by the centralizer of the map $v\mapsto D\Phi(0)v$ in the space of non-linear formal maps.
This generalizes the result of Wilhelm Klingenberg Jr.(see above) ($p=1$) in the sense that, there are holomorphic coordinates such that $(M_Q)$ intersects $\{z_2=0\}$ along $\{(\zeta_1,\eta_1)\in \Bbb R^2\;|\;\zeta_1 \eta_1=0\}$.

The core of these problems rests on geometric properties of dynamical systems
associated to each situation. To be more
specific, we shall deal, in the first part of this article, with germs of holomorphic diffeomorphisms of
$(\Bbb C^n,0)$ in a neighbourhood of the origin (a common fixed point). We
shall consider those whose linear part at the origin is different from the
identity. The main point is a result which gives the existence of germ of analytic subset of
$(\Bbb C^n,0)$ invariant by a abelian group of such diffeomorphisms under some
diophantine condition. This kind of result was obtained by the author for a
germ of holomorphic vector field at singular point \cite{Stolo-dulac}.

%%% Local Variables: 
%%% mode: latex
%%% TeX-master: "webster"
%%% TeX-master: "webster"
%%% TeX-master: "webster"
%%% TeX-master: "webster"
%%% TeX-master: "webster"
%%% End: 

\section{Abelian group of diffeomorphisms of $(\Bbb C^n,0)$ and their invariant sets}

The aim of this section is to prove the existence of complex analytic invariant subset for a commuting family of germs of holomorphic diffeomorphisms in a neighbourhood of a common fixed point. This is very inspired by a previous article of the author concerning holomorphic vector fields. Althought the objects are not the same, some of the computations are identical and we shall refer to them when possible.

Let $D_1:=\text{diag}(\mu_{1,1},\ldots,\mu_{1,n}),\ldots, D_l:=\text{diag}(\mu_{l,1},\ldots,\mu_{l,n})$ be diagonal invertible matrices.
Let us consider a familly $F:=\{F_i\}_{i=1,\ldots l}$ of {\bf commuting germs of holomorphic diffeomorphisms} of $(\Bbb C^n,0)$ which linear part, at the origin, is $D:=\{D_ix\}_{i=1,\ldots l}$ : 
$$
F_i(x)=D_ix+f_i(x), \quad \text{with}\quad f_i(0)=0,\;Df_i(0)=0,\;f_i\in {\cal O}_n.
$$

Let ${\cal I}$ be an ideal of ${\cal O}_n$ {\bf generated by monomials} of
$\Bbb C^n$. Let $V({\cal I})$ be the germ at the origin, of the analytic
subset of $(\Bbb C^n,0)$ defined by ${\cal I}$. It is left invariant by
the familly $D$. Let us set $\hat {\cal I}:=\widehat{\cal O}_n\otimes {\cal
  I}$. Here we denote ${\cal O}_n$ (resp. $\widehat{\cal O}_n$) the ring of
germ of holomorphic function at the origin (resp. ring of formal power series)
of $\Bbb C^n$. Let $Q=(q_1,\ldots, q_n)\in \Bbb N^n$ and
$x=(x_1,\ldots,x_n)\in \Bbb C^n$, we shall write
$$
|Q|:=q_1+\cdots +q_n,\quad x^Q:=x_1^{q_1}\cdots x_n^{q_n}.
$$ 

Let $\{\omega_k(D,{\cal I})\}_{k\geq 1}$ be the sequence of positive numbers defined by
$$
\omega_k(D,{\cal I})=\inf\left\{\max_{1\leq i\leq l}|\mu_{i}^Q-\mu_{i,j}|\neq 0\;|\;2\leq |Q|\leq 2^k, 1\leq j\leq n,Q\in \Bbb N^n, x^Q\not\in {\cal I}\right\}.
$$
Let $\{\omega_k(D)\}_{k\geq 1}$ be the sequence of positive numbers defined by
$$
\omega_k(D,{\cal I})=\inf\left\{\max_{1\leq i\leq l}|\mu_{i}^Q-\mu_{i,j}|\neq 0\;|\;2\leq |Q|\leq 2^k, 1\leq j\leq n,Q\in \Bbb N^n \right\}.
$$
\begin{defis}
\begin{enumerate}
\item We shall say that the ideal ${\cal I}$ is properly embedded if it has a set of monomial generators not involving a nonempty set ${\cal S}$ of variables.
\item We shall say that the familly $D$ is {\bf diophantine} (resp. on ${\cal I}$) if
$$
-\sum_{k\geq 1}\frac{\ln \omega_k(D)}{2^k}<+\infty \quad (resp. -\sum_{k\geq 1}\frac{\ln \omega_k(D,{\cal I})}{2^k}<+\infty).
$$
\item We shall say that the familly $F$ is {\bf formally linearizable on
    $\hat{\cal I}$} if there exists a formal diffemorphism $\hat\Phi$ of
  $(\Bbb C^n,0)$, tangent to the identity at the origin such that
  $\hat\Phi_*F_i-D_ix \in (\hat{\cal I})^n$ for all $1\leq i\leq l$.
\item A linear anti-holomorphic involution of $\Bbb C^n$ is a map
$\rho(z)=P\bar z$ where the matrix $P$ satisfies $P\bar P=Id$; $\bar z$ denotes the complex conjugate of $z$.
\item We shall say that ${\cal I}$ is {\bf compatible} with a anti-linear involution $\rho$ if the map $\rho^*:\widehat{\cal O}_{n+p-1}\rightarrow \overline{\widehat{\cal O}_{n+p-1}}$) defined $\rho^*(f)=f\circ \rho$  maps  $\widehat{\cal I}$ to $\overline{\widehat{\cal I}}$ and  $\widehat{CI}$ to $\overline{\widehat{CI}}$.
\end{enumerate}
\end{defis}
Let $\widehat{\cal O}_n^D$ be the ring of formal invariant of the familly $D$, that is 
$$
\widehat{\cal O}_n^D:=\{f\in \widehat{\cal O}_n\,|\;f(D_ix)=f(x)\; i=1,\ldots, l  \}.
$$
It can be shown (as in proposition 5.3.2 of \cite{stolo-ihes}) that this ring
is generated by a finite number of monomials $x^{R_1},\ldots, x^{R_p}$ and
that the non-linear centralizer ${\cal C}_D$ of $D$ is a module over
$\widehat{\cal O}_n^D$ of finite type. Let $\text{ResIdeal}$ be the ideal
generated by the monomials $x^{R_1},\ldots, x^{R_p}$ in ${\cal O}_n$.

\begin{theos}\label{theo-invariant}
Let ${\cal I}$ be a monomial ideal (resp. properly embedded). Assume that the familly $D$ is diophantine (resp. on
${\cal I}$). If the familly $F$ is formally linearizable on $\hat{\cal
  I}$, then it is holomorphically linearizable on ${\cal I}$. Moreover,
there exists a unique such a diffeomorphism $\Phi$ such that the projection of the Taylor expansion of $\Phi-Id$ onto ${\cal I}\cup {\cal C}_D$ vanishes.

Moreover, let $\rho$ be a linear anti-holomorphic involution such that $\rho {\cal C}_D\rho={\cal C}_D$. We assume that  ${\cal I}$ is compatible with $\rho$. Assume that, for all $1\leq i\leq l$, $\rho\circ F_i\circ \rho$ belongs to the group generated by the $F_i$'s. Then $\Phi$ and $\rho$ commute with each other.
\end{theos}
This theorem can be rephrased as follow~: Under the afore-mentioned
diophantine condition,  then there exists a germ of holomorphic diffeomorphism
$\Phi$ such that $\Phi_*F_i-D_ix \in ({\cal I})^n$ for all $1\leq i\leq l$. As
a consequence, in a good holomorphic coordinates system, the analytic subset
$V({\cal I})$ is left invariant by each $F_i$ and its restriction to it is the
linear mapping $x\mapsto D_{i|V({\cal I})}x$. 
\begin{rems}
The familly $D$ can be diophantine while none of the $D_i$'s is.
\end{rems}
The second part of the theorem will be used for applications in the third part of the article.
\begin{coros}\label{invar-facile}
\begin{enumerate}
\item If the ring of invariant of $D$ reduces to the constants and if $D$ is
  diophantine, then $F$ is holomorphically linearizable in a neighbourhood of
  the origin. For one diffeomorphism, this was obtain by H. R\"ussmann
  \cite{russmann-ihes,russmann-diffeo} and by T. Gramtchev and M.Yoshino
  \cite{gram-yoshi} for an abelian group under a slightly coarser diophantine condition.
\item The existence of an invariant manifold for a germ of diffeomorphism was obtain by J. P\"oschel \cite{Posch}. Despite the fact that we are dealing with a family of diffeomorphisms, the main difference is that we are able to linearize {\bf simultaneously} on each irreducible component of analytic set.
\end{enumerate}
\end{coros}
According to M. Chaperon \cite{chaperon-ast}[theorem 4, p.132], if the family of diffeomorphisms is abelian then there exists a formal diffeomorphism $\hat \Phi$ such that
$$
\hat \Phi_*F_i(D_jz)= D_j\hat \Phi_*F_i(z),\quad 1\leq i,j\leq l.
$$
We call the family of $\hat \Phi_*F_i$'s a {\bf formal normal form} of the family $F$.
Then we have the following corollary~:
\begin{coros}\label{invar-group}
Let $F$ be an abelian family of germs of holomorphic diffeomorphisms of $(\Bbb C^n,0)$.  Let us assume that $D$ is diophantine on $ResIdeal$. If the non-linear centralizer of $D$ is generated by the $x^{R_i}$'s then $F$ is holomorphically linearizable on ${\cal I}$.
\end{coros}
\begin{rems}
The condition that the non-linear centralizer of $D$ is generated by the $x^{R_i}$'s means: if $\mu_i^Q=\mu_{i,j}$ for some $Q\in \Bbb N^n_2$, $1\leq j\leq n$ and for all $1\leq i\leq l$, then 
$x^Q$ belongs to the ideal generated by $x^{R_1},\ldots, x^{R_p}$. This a very weak condition since only all but a finite number of resonances satisfy this condition.
\end{rems}
%Let us consider the ring of invariant of $D$ :
%$$
%\widehat{\cal O}_n^D:=\{h\in\widehat{\cal O}_n\;|\;h\circ D_i(x)=h(x),\;i=1,\ldots,l\}.
%$$
%If not reduced to the constant, this ring is generated as an algebra by a finite number of monomials. This %mean that there exists monomials $x^{R_1},\ldots, x^{R_p}$ such that
%$$
%\widehat{\cal O}_n^D=\Bbb C[[x^{R_1},\ldots, x^{R_p}]].
%$$
%Let us assume that the centralizer of $D$ in the subspace of elements of $\widehat{\cal O}_n^n$ of order %greater or equal than $2$ is generated by $x_i\frac{\partial}{\partial x_i}$ over $\widehat{\cal O}_n^D$. This %means that if $\mu_i^Q=\mu_{i,j}$, for all $1\leq i\leq n$, where $Q\in \Bbb N^n_2$ and $1\leq j \leq n$ , %then $Q=l_1R_1+\cdots l_pR_p+E_j$ for some non-negative integers $l_i$ ($E_j$ denotes the $j$th vector of the %natural basis of $\Bbb Z^n$).
%Let ${\cal I}$ be the ideal of ${\cal O}_n$ generated by the monomials $x^{R_1},\ldots, x^{R_p}$. We shall %denote by $V({\cal I})$ the germ at the origin of the analytic set defined by ${\cal I}$.

We shall prove that there exists a holomorphic map $\phi : (\Bbb C^n,0)\rightarrow (\Bbb C^n,0)$, tangent to the identity at the origin, such that
$$
\Phi^{-1}\circ F_i \circ \Phi (y)= G_i(y):= D_iy +g_i(y)\quad i=1,\ldots, l
$$ 
where the components of $g_i$ are non-linear holomorphic functions and belong
to the ideal ${\cal I}$. It is unique if we require that its projection on
${\cal I}\cup \text{ResIdeal}$ is zero.

Let us set $x_j=\Phi_j(y):= y_j+\phi_j(y)$, $j=1,\ldots, n$. Let us expand the equations $F_i \circ \Phi (y)=\Phi\circ G_i$, $i=1,\ldots,l$.
For all $1\leq j\leq n$ and all $i=1,\ldots,l$, we have
\begin{eqnarray*}
\mu_{i,j}y_j+g_{i,j}(y)+\phi_{j}(G_i(y)) & = & \mu_{i,j}(y_j+\phi_{j}(y))+f_{i,j}(\Phi (y))\\
g_{i,j}(y)+\phi_{j}(D_iy) & = & \mu_{i,j}\phi_{j}(y)+f_{i,j}(\Phi (y))\\
+(\phi_{j}(G_i(y))-\phi_{j}(D_iy))& & \\
\end{eqnarray*}
Let us expand the functions at the origin~:
$$
f_{i,j}(y)=\sum_{Q\in\Bbb N^n_2}f_{i,j,Q}y^Q,\;g_{i,j}(y)=\sum_{Q\in\Bbb N^n_2}g_{i,j,Q}y^Q\text{ and }\phi_{j}(y)=\sum_{Q\in\Bbb N^n_2}\phi_{j,Q}y^Q.
$$
Then we have
\begin{equation}
\sum_{Q\in\Bbb N^n_2}\delta^i_{Q,j}\phi_{j,Q}y^Q + g_{i,j}(y) = f_{i,j}(\Phi (y))-(\phi_{j}(G_i(y))-\phi_{j}(D_iy))\label{conjugaison}
\end{equation}
where 
$$
\delta^i_{Q,j} := \mu_i^Q-\mu_{i,j},\quad\mu_i :=(\mu_{i,1},\ldots,\mu_{i,n}).
$$
%We emphasize that if $y^Q$ does not belong to the ideal ${\cal I}\subset ResIdeal$ then, for any $1\leq j\leq n$, there exists $1\leq i_0\leq n$ such that $\delta^{i_0}_{Q,j} \neq 0$.

Let $\{f\}_Q$ denotes the coefficient of $x^Q$ in the Taylor expansion at the origin of $f$.
We compute $\phi_{j,Q}$ and $g_{i,j,Q}$ by induction on $|Q|\geq 2$ in the following way :
\begin{itemize}
\item if $y^Q$ does not belongs to ${\cal I}$ and $\max_i |\delta^{i}_{Q,j}|\neq 0$, then there exists $1\leq i_0\leq l$ such that $|\delta^{i_0}_{Q,j}|=\max_i |\delta^{i}_{Q,j}|$. We set 
\begin{eqnarray*}
\phi_{j,Q} & = & \frac{1}{\delta^{i_0}_{Q,j}}\left\{f_{i_0,j}(\Phi (y))-(\phi_{j}(G_{i_0}(y))-\phi_{j}(D_{i_0}y))\right\}_Q\\
g_{i,j,Q} & = & 0.
\end{eqnarray*}
\item If $y^Q$ does not belongs to ${\cal I}$ and $\max_i |\delta^{i}_{Q,j}|=0$, then we have 
$$
\left\{f_{i_0,j}(\Phi (y))-(\phi_{j}(G_{i_0}(y))-\phi_{j}(D_{i_0}y))\right\}_Q=0
$$
and we set $\phi_{j,Q} = 0 = g_{i,j,Q}$.

\item If $y^Q$ belongs to ${\cal I}$, we set 
\begin{eqnarray*}
\phi_{j,Q} & = & 0\\
g_{i,j,Q} & = & \left\{f_{i,j}(\Phi (y))-(\phi_{j}(G_{i}(y))-\phi_{j}(D_{i}y))\right\}_Q.
\end{eqnarray*}
\end{itemize}
\begin{lemms}
The formal diffeomorphism $\Phi$ defined above linearizes simultanueously the family $F$ on $\hat {\cal I}$ where $\hat {\cal I}:=\widehat{\cal O}_n\otimes {\cal I}$.
\end{lemms}
\begin{proof}
For all $1\leq i,j\leq j\leq l$, we have 
$$
F_i\circ F_j  =  F_j\circ F_i\;\;\text{ thus }\;\; F_i\circ F_j\circ \Phi  =  F_j\circ F_i\circ \Phi.
$$
Therefore, we have
$$
D_iD_j\Phi+D_i(f_j\circ \Phi) +f_i(D_j\Phi+f_j\circ \Phi) =  D_jD_i\Phi+D_i(f_i\circ \Phi) +f_j(D_i\Phi+f_i\circ \Phi),
$$
so that
\begin{eqnarray*}
f_j(\Phi\circ D_i)-D_i(f_j\circ \Phi) & = & f_i(\Phi\circ D_j)-D_j(f_i\circ \Phi)\\
+ f_j(D_i\Phi+f_i\circ \Phi)-f_j(\Phi\circ D_i)& & + f_i(D_j\Phi+f_j\circ \Phi)-f_i(\Phi\circ D_j).
\end{eqnarray*}
Moreover, we have $F_i\circ \Phi = \Phi\circ G_i$. Hence, we have
\begin{eqnarray*}
f_i(D_j\Phi+f_j\circ \Phi)-f_i(\Phi\circ D_j) & = & f_i\circ F_j\circ \Phi -f_i(\Phi\circ D_j) \\
&= & f_i\circ \Phi\circ G_j-f_i\circ\Phi\circ D_j \\
& = & D(f_i\circ \Phi)(D_jy)g_j+\cdots\\
\end{eqnarray*}
Assume the $F_i$'s are linearized on $V({\cal I})$ up to order $k\geq 2$. This means that, for any $1\leq m\leq n$ and any $1\leq i\leq l$, the $k$-jet $J^k(g_{i,m})$ belongs to ${\cal I}$. The previous computation shows that the $k+1$-jet of $f_i(D_j\Phi+f_j\circ \Phi)-f_i(\Phi\circ D_j)$ depends only on the $k$-jet of $g_j$ and belongs to ${\cal I}$. The same is true for $\phi_{j}(G_{i}(y))-\phi_{j}(D_{i}y)$. Therefore, if $Q\in \Bbb N^n_2$ with $|Q|=k+1$ is such that $x^Q$ does not belong to ${\cal I}$, then we have 
$$
\{f_j(\Phi\circ D_i)-D_i(f_j\circ \Phi)\}_Q  =  \{f_i(\Phi\circ D_j)-D_j(f_i\circ \Phi)\}_Q;
$$
that is, for all $1\leq m\leq n$,
$$
(\mu_{i}^Q-\mu_{i,m})\{f_{j,m}\circ \Phi\}_Q  = (\mu_{j}^Q-\mu_{j,m})\{f_{i,m}\circ \Phi\}_Q.
$$
This means that equation $(\ref{conjugaison})$ is solved by induction and that $\Phi$ linearizes formally the $F_i$'s on $V(\hat{\cal I})$.
\end{proof}

Let $\rho$ be a linear anti-holomorphic involution satisfying the assumptions of the
theorem. We have $F_i\circ\Phi = \Phi\circ G_i$ where $G_i$ is linearized along $\widehat{\cal I}$. Hence, we have 
$$
(\rho \circ F_i\circ \rho)\circ(\rho\circ\Phi\circ\rho) = (\rho\circ\Phi\circ\rho)\circ (\rho \circ G_i\circ \rho).
$$
Let us set $\tilde F_i:=\rho \circ F_i\circ \rho$. By assusmptions, $\tilde F_i$ belongs to the group generated by the $F_i$'s. Since $\rho^*\widehat{\cal I}\subset \overline{\widehat{\cal I}}$, then $\rho \circ G_i\circ \rho$ is a formal
diffeomorphism which is linearized on $\widehat{\cal I}$. 
%Its restriction to it is  equal to the mapping $x\mapsto D\tilde F_i(0)x$. 
By assumptions, the projection of $\rho\circ\Phi\circ\rho-Id $ onto ${\cal I}\cup {\cal C}_D$ vanishes
identically.  By uniqueness, we have $\rho\circ\Phi\circ\rho=\Phi$ since $\Phi$ linearizes $\tilde F_i$ on ${\cal I}$.

We shall prove, by using the majorant method, that $\Phi$ actually converges
on a polydisc of positive radius centered at the origin. Let us define $\Bbb
N_2^n\setminus \hat{\cal I}$ to be the set of multiindices $Q\in \Bbb N^n$
such that $|Q|\geq 2$ and $x^Q\not\in \widehat{\cal I}$. Let $f=\sum_Q f_Qx^Q$
and $g=\sum_Q g_Qx^Q$ be formal power series. We shall say that $g$ {\it
  dominates} if $|f_Q|\leq |g_Q|$ for all multiindices $Q$.

First of all, for all $1\leq j\leq n$ and all $ Q\in \Bbb N_2^n\setminus \hat{\cal I}$ such that $\max_{1\leq i\leq l}|\delta_{j,Q}^i|\neq 0$, we have
$$
|\phi_{j,Q}||\delta_{j,Q}|=|\{f_{i_0(Q),j}(\Phi)\}_Q|\leq \{\bar f_{i_0(Q),j}(y+\bar \phi)\}_Q
$$
where $|\delta_{j,Q}|=\max_{1\leq i\leq l}|\delta_{j,Q}^i|=|\delta_{i_0(Q,j),j,Q}|$. In fact, $ \{f_i\circ \Phi\circ G_j-f_i\circ\Phi\circ D_j \}_Q=0$ whenever $ Q\in \Bbb N_2^n\setminus \hat{\cal I}$. This inequality still holds  if $\max_{1\leq i\leq l}|\delta_{j,Q}^i|= 0$. Let us set 
\begin{itemize}
\item $\delta_{Q}:= \min\{|\delta_{j,Q}|, {1\leq j\leq n}$ such that  $\delta_{j,Q}\neq 0\}$,\\
\item $\delta_{Q}:=0$ if $\max_{1\leq i\leq l}|\delta_{j,Q}^i|= 0$.
\end{itemize}

Let us sum over $1\leq j\leq n$ the previous inequalities. We obtain for all $Q\in \Bbb N_2^n\setminus \hat{\cal I}$,
$$
\delta_Q\sum_{j=1}^n{|\phi_{j,Q}|}\leq\sum_{j=1}^n{|\phi_{j,Q}||\delta_{j,Q}|}\leq \left\{\sum_{j=1}^n\bar f_{i_0(Q,j),j}(y+\bar \phi)\right\}_Q\leq \left\{\sum_{i=1}^l\left(\sum_{j=1}^n\bar f_{i,j}\right)(y+\bar \phi)\right\}_Q.
$$

Since $\sum_{i=1}^l\sum_{j=1}^n f_{i,j}$ vanishes at the origin with its derivative as well, there exists positives constants $a,b$ such that
$$
\sum_{i=1}^l\sum_{j=1}^n{f_{i,j}}\prec \frac{a\left(\sum_{j=1}^n{x_j}\right)^2}{1-b\left(\sum_{j=1}^n{x_j}\right)}.
$$
Since the Taylor expansion of the right hand side has non-negative coefficients, we obtain
$$
\delta_Q\tilde\phi_Q\leq \left\{\frac{a\left(\sum_{j=1}^n{y_j+\tilde \phi}\right)^2}{1-b\left(\sum_{j=1}^n{y_j+ \tilde \phi}\right)}\right\}_Q
$$
where we have set $\tilde\phi_Q:=\sum_{j=1}^n{|\phi_{j,Q}|}$ and $\tilde \phi=\sum_{Q\in \Bbb N_2^n}\tilde\phi_Qx^Q$. Here, we have set $\tilde\phi_Q=0$ whenever $\delta_Q=0$. 

Let us define the formal power series $\sigma(y)=\sum_{Q\in \Bbb N_2^n}{\sigma_Qy^Q}$ as follow :
\begin{eqnarray*}
\forall Q\in \Bbb N_2^n\setminus (\Bbb N_2^n\setminus \hat{\cal I})\;\;\; \sigma_Q & = & 0\\
\forall Q\in \Bbb N_2^n\setminus \hat{\cal I}\;\;\; \sigma_Q & = & \left\{\frac{a\left(\sum_{j=1}^n{y_j+ \sigma}\right)^2}{1-b\left(\sum_{j=1}^n{y_j+ \sigma}\right)}\right\}_Q
\end{eqnarray*}
\begin{lemms}\cite{Stolo-dulac}[Lemme 2.1]
The series $\sigma$ is convergent in a neighbourhood of the origin $0\in \Bbb C^n$.
\end{lemms}

Let us define the sequence $\{\eta_Q\}_{Q\in \Bbb N^n_1\setminus \hat{\cal I}}$ of positive number  as follow~:\\
\ben
\item $\forall P\in \Bbb N^n_1\setminus \hat{\cal I}$ tel que $|P|=1$, $\eta_{P}  =  1$ ( such multiindice exists. ),
\item $\forall Q\in \Bbb N^n_2\setminus \hat{\cal I}$ with $\delta_Q\neq 0$
$$
\delta_Q \eta_Q  = \max_{{\substack{ Q_j\in \Bbb N^n_1,S\in \Bbb N^n  \\
Q_1+\cdots+Q_p+S=Q  }}}{\eta_{Q_1}\cdots\eta_{Q_p}},
$$
the maximum been taken over the sets of $p+1$, $1\leq p\leq |Q|$, multiindices $Q_1,\ldots,Q_p,S$ such that
$\forall 1\leq j\leq p,\;Q_j\in \Bbb N_1^n,\; |Q_j|<|Q|$, $S\in \Bbb N^n$. These sets are not empty.
\item $\forall Q\in \Bbb N^n_2\setminus \hat{\cal I}$ with $\delta_Q= 0$, $\eta_Q =0$.
\een
This sequence is well defined. In fact, if $Q\in \Bbb N^n_2\setminus \hat{\cal I}$, then there exists
 multiindices $Q_1,\ldots,Q_p,S$ such that $Q=Q_1+\ldots+Q_p+S$, $\forall 1\leq j\leq p,\;Q_j\in \Bbb N_1^n,\; |Q_j|<|Q|,\;S\in \Bbb N^n$.
In this case, $\forall 1\leq j\leq p,\;\;Q_j\in \Bbb N_1^n\setminus \hat{\cal I}$.

The following lemmas are the key points.
\begin{lemms}\cite{Stolo-dulac}[Lemme 2.2]\label{majoration}
For all $Q\in \Bbb N^n_2\setminus \hat{\cal I}$, we have  $\tilde \phi_Q\leq \sigma_Q\eta_Q$.
\end{lemms}
\begin{lemms}\cite{Stolo-dulac}[Lemme 2.3]\label{lemm-divis}
There exists a constant $c>0$ such that $\forall Q\in \Bbb N_2^n\setminus \hat{\cal I},\;\;\eta_Q\leq c^{|Q|}$.
\end{lemms}

Let $\theta>0$ be such that $4\theta:=\min_{i,j}|\lambda_{i,j}|\leq 1$ (we can
always assume this, even if this means using the inverse of one of the  diffeomorphisms). If the ideal ${\cal I}$ is properly embedded, then we shall set 
$$
4\theta:=\min_{1\leq i\leq l, j\in{\cal S}}|\lambda_{i,j}|\leq 1
$$ 
where ${\cal S}$ denotes the set of variables not involves in any generator.
In particular, we have the property that if $x^Q\not\in {\cal I}$ then $x_sx^Q\not\in {\cal I}$ for all $s\in {\cal S}$.

By definition,  $\eta_Q$ is a product of $1/\delta_{Q'}$ with $|Q'|\leq|Q|$. Let $k$ be a non-negative integer. Let us define $\phi^{(k)}(Q)$ (resp. $\phi^{(k)}_j(Q)$) to be the number of $1/\delta_{Q'}$'s present in this product and such that $0\neq \delta_{Q'}<\theta\omega_k(D,{\cal I})$ (resp. and $\delta_Q=\delta_{j, Q}$). The lemma is a consequence of the following proposition
\begin{prop}\cite{Stolo-dulac}[lemme 2.8]\label{estim-bruno}
For all $Q\in {\Bbb N}^n_2\setminus \hat{\cal I}$, we have $\phi^{(k)}(Q)\leq 2n\frac{|Q|}{2^k}$ if $|Q|\geq 2^k+1$; and $\phi^{(k)}(Q)=0$ if $|Q|\leq 2^k$.
\end{prop}

In fact, $\phi^{(k)}(Q)$ bounds the number of $1/\delta_{Q'}$'s appearing in the product defining $\eta_Q$ and such that $\theta\omega_{k+1}(D, {\cal I})\leq \delta_{Q'}<\theta\omega_k(D, {\cal I})$.
\begin{proof}[Proof of lemma \ref{lemm-divis}]
Let $r$ be the integer such that $2^r+1\leq |Q|<2^{r+1}+1$. Then we have
$$
\eta_Q\leq \prod_{k=0}^r{\left(\frac{1}{\theta\omega_{k+1}(D, {\cal I})}\right)^{\phi^{(k)}(Q)}}.
$$
By applying the Logarithm and proposition \ref{estim-bruno}, we obtain
\begin{eqnarray*}
\ln \eta_Q &\leq &\sum_{k=0}^l{2n\frac{|Q|}{2^k}\left(\ln \frac{1}{\theta\omega_{k+1}(D)}\right)}\\
& \leq & |Q|\left(-2n\sum_{k\geq 0}{\frac{\ln \omega_{k+1}(D)}{2^k}}+2n\ln \theta^{-1}\sum_{k\geq 0}{\frac{1}{2^k}}\right).
\end{eqnarray*}
Since the familly $D$ is diophantine, we obtain 
$\eta_Q \leq  c^{|Q|}$ for some positive constant $c$.
\end{proof}

For any positive integer $k$, for any $1\leq j\leq n$, let us consider the function defined on ${\Bbb N}^n_2\setminus \hat{\cal I}$ to be
$$
\forall Q\in \Bbb N_2^n\setminus \hat{\cal I},\;\;\;
\psi^{(k)}_j(Q)=\left \{ \begin{array}{l} 1 \;\;\;\;\mbox{ if }\; \delta_Q=|\delta_{j,Q}|\neq 0 \;\mbox{ and }\; |\delta_{j,Q}|<\theta\omega_k(D, {\cal I})\\
0 \;\;\;\;\mbox{ if }\; \delta_Q=0 \;\mbox{ or }\;\delta_Q\neq |\delta_{j,Q}|\;\mbox{ or }\; |\delta_{j,Q}|\geq \theta\omega_k(D, {\cal I})\\ \end{array} \right .
$$
Then we have, 
$$
0\leq \phi_j^{(k)}(Q)  = \psi^{(k)}_j(Q)+\max_{{\substack{ Q_j\in \Bbb N^n_1,S\in \Bbb N^n  \\
Q_1+\cdots+Q_p+S=Q  }}}{\left(\phi_j^{(k)}(Q_1)+\cdots+\phi_j^{(k)}(Q_p)\right)}.
$$
The proof of propositon \ref{estim-bruno} identitical to the proof of \cite{Stolo-dulac}[lemme 2.8] except that we have to use the following version of \cite{Stolo-dulac}[lemme 2.7].

\begin{lemms}\label{lemm-tech}
Let $Q\in {\Bbb N}^n_2\setminus \hat{\cal I}$ be such that $\psi^{(k)}_j(Q)=1$.
If $Q=P+P'$ with $(P,P')\in \Bbb N^n_1\times\Bbb N^n_2$ and $|P|\leq 2^k-1$, 
then $(P,P')\in \Bbb N^n_1\setminus \hat{\cal I}\times\Bbb N^n_2\setminus \hat{\cal I}$ and $\psi^{(k)}_j(P')=0$.
\end{lemms}

\begin{proof}  Clearly, if $Q=P+P'\in {\Bbb N}^n_2\setminus \hat{\cal I}$ then $(P,P')\in \Bbb N^n_1\setminus \hat{\cal I}\times\Bbb N^n_2\setminus \hat{\cal I}$.
There are two cases to consider :
\ben
\item if $\delta_{P'}\neq |\delta_{j,P'}|$ or $\delta_{P'}=0$ then $\psi^{(k)}_j(P')=0$, by definition.
\item if $\delta_{P'}=|\delta_{j,P'}|\neq 0$, assume
that $\delta_{P'}<\theta\omega_k(D, {\cal I})$. Then, for all $1\leq i\leq l$, we have
$$
|\lambda_i^{P'}|>|\lambda_{i,j}|-\theta\omega_k(D, {\cal I})\geq 4\theta -2\theta= 2\theta.
$$
It follows that, for all $1\leq i\leq l$,
\begin{eqnarray*}
2\theta\omega_k(D, {\cal I}) & > & |\lambda_i^Q-\lambda_{i,j}|+|\lambda_i^{P'}-\lambda_{i,j}|\\
& > & |\lambda_i^Q-\lambda_i^{P'}| = |\lambda_i^{P'}||\lambda_i^{P}-1|.
\end{eqnarray*}
If ${\cal I}$ is properly embedded, for all $a\in {\cal S}$, we have $x_ax^P\not\in {\cal I}$. 
Therefore,  for all $1\leq i\leq l$, we have
\begin{eqnarray*}
2\theta\omega_k(D, {\cal I}) & > & 2\theta|\lambda_{i,a}|^{-1}|\lambda_i^{P+E_a}-\lambda_{i,a}|\\
& > & 2\theta|\lambda_{i,a}|^{-1}\omega_k(D, {\cal I}).
\end{eqnarray*}
This contradicts the facts that $\min_{1\leq i\leq l,a\in {\cal S}}|\lambda_{i,a}|\leq 1$. 
If ${\cal I}$ is not properly embedded, then we obtain 
$$
2\theta\omega_k(D) > 2\theta|\lambda_{i,a}|^{-1}\omega_k(D)
$$
for all $1\leq a\leq n$. It is still a contradiction.
\een
Hence, we have shown that $\psi^{(k)}_j(P')=0$. 
\end{proof}

%%% Local Variables: 
%%% mode: latex
%%% TeX-master: "webster"
%%% End: 

\section{Family of totally real $n$-manifolds in $(\Bbb C^n,0)$} 

Let us consider a family $M:=\{M_i\}_{i=1,\ldots, m}$ of real analytic totally real $n$-submanifold of $\Bbb C^n$ passing throught the origin. Locally, each $M_i$ is the fixed point set of an anti-holomorphic involution $\rho_i$~: $M_i = FP(\rho_i)$ and $\rho_i\circ \rho_i=Id$. This means that
$$
\rho_i(z) := B_i\bar z + R_i(\bar z)
$$
where $R_i$ is a germ of holomorphic function at the origin with $R_i(0)=0$ and $DR_i(0)=0$. Each matrix $B_i$ is invertible and satisfies $B_i\bar B_i=Id$.
The tangent space, at the origin, of $M_i$ is the totally real $n$-plane
$$
\{z=B_i\bar z\}
$$
We assume that there are all distinct one from another.
Their intersection at the origin is the set
$$
\left\{z\in \Bbb C^n\;|\;B_i\bar z=z,\;i=1,\ldots,m\right\}\subset \left\{z\in \Bbb C^n\;|\;B_i\bar B_j z=z,\;i,j=1,\ldots,m\right\}.
$$
It is contained in the common eigenspace of the $B_i\bar B_j$'s associated to the eigenvalue $1$. {\bf We shall not assume that this space is reduced to $0$}.

Let us consider the group $G$ generated by the germs of holomorphic diffeomorphisms of $(\Bbb C^n,0)$
$F_{i,j}:=\rho_i\circ \rho_j$, $1\leq i,j\leq m$. Let $D_{i,j}:=B_i\bar B_j$ be the linear part at the origin of $F_{i,j}$. Let us set 
$$
F_{i,j}:=D_{i,j}z+f_{i,j}(z)
$$
where $f_{i,j}$ is a germ of holomorphic function at the origin with $f_{i,j}(0)=0$ and $Df_{i,j}(0)=0$.

Let us write the relation $F_{i,j}=\rho_i\circ \rho_j$ and $\rho_i\circ \rho_i=Id$. We obtain
\begin{eqnarray}
f_{i,j}(z) & = & B_i\bar R_j(z) + R_i(\bar \rho_j)\label{equ1}\\
0 & = & B_i\bar R_i(z) + R_i(\bar \rho_i).\label{equ2}
\end{eqnarray}
By multiplying the first equation by $\bar B_i$, we obtain
$$
\bar R_j(z) = \bar B_i f_{i,j}(z)- \bar B_i R_i(\bar \rho_j).
$$
Hence,we have
$$
0  =   B_j\bar B_i f_{i,j}(z)- B_j\bar B_i R_i(\bar \rho_j)+  B_i \bar f_{i,j}(\bar \rho_j)-  B_i \bar R_i(\rho_j\circ\rho_j).
$$
Let us mupltiply by $\bar B_i$ on the left and take the conjugation. We obtain
$$
0= D_{i,j}B_i\bar f_{i,j}(\bar z)-D_{i,j}B_i \bar R_i(\rho_j)+ f_{i,j}(\rho_j)-R_i(\bar z).
$$
On the other hand, by evaluating equation $(\ref{equ2})$ at $\bar \rho_j$, we obtain
$$
0=B_i\bar R_i(\bar \rho_j) + R_i(\bar F_{i,j}).
$$
At the end,we obtain
\begin{equation}
R_i(\bar z) - D_{i,j}R_i(\bar F_{i,j}) =  D_{i,j}B_i\bar f_{i,j}(\bar z)+f_{i,j}(\rho_j)\label{conj}.
\end{equation}
\begin{defis}
The $\rho_i$'s are simultaneously normalizable whenever $R_i(\bar z) - D_{i,j}R_i(\bar D_{i,j}\bar z)=0$ for all $1\leq i,j\leq l$.
\end{defis}
\begin{rems}\label{rem-normalisation}
If the group $G$ is holomorphically linearizable at the origin then the $\rho_i$'s are simultaneously normalizable.

Moreover, asssume the $D_{i,j}$'s are simultaneously diagonalizable and let us set $D_{i,j}=\text{diag}(\mu_{i,j,k})$. Then, for any $1\leq k\leq n$ and any $1\leq j\leq m$, the $k$-component $\rho_{i,k}$ of $\rho_i$ can be written as
$$
\left(\rho_{i}(z)-B_i\bar z\right)_k = \sum_{\substack{Q\in \Bbb N^n_2\\ \forall j,\;\bar\mu_{i,j}^Q=\mu_{i,j,k}^{-1}}} \rho_{i,k,Q}\bar z^Q.
$$
Here, $(f)_k$ denotes the kth-component of $f$.
\end{rems}
As a consequence, we have
\begin{theos}
Let us assume that the group $G$ associated to the family of totally real
submanifolds $M$ is a semi-simple Lie group. Then the $\rho_i$'s are
simultaneously and holomorphically normalizable in a neighbourhood of the origin.
\end{theos}
\begin{proof}
It is classical \cite{Kushnirenko, stern-semi-simple, Ghys} that if the Lie
group $G$ of germs of diff\'eomorphisms at a common fixed point is semi-simple then it is holomorphically linearizable in a neighbourhood of the origin. Then, apply the previous remark \ref{rem-normalisation}.
\end{proof}
\begin{defis}
We shall say that such a family $M=\{M_i\}_{i=1,\ldots, m}$ of totally real $n$-submanifold of $(\Bbb C^n,0)$ intersecting at the origin is {\bf commutative} if the group $G$ is abelian.
\end{defis}
>From now on, we shall assume that {\bf $M$ is commutative} and that the {\bf family $D$}
of linear part of the group $G$ at the origin {\bf is diagonal}. In other words, $D_{i,j}=\text{diag}(\mu_{i,j,k})$.
Let ${\cal I}$ be a monomial ideal of ${\cal O}_n$. It is genrated by some
monomials $x^{R_1},\ldots, x^{R_p}$. We shall denote $\bar{\cal I}$ the ideal
of $\Bbb C[[\bar x_1,\ldots,\bar x_n]]$ generated by $\bar x^{R_1},\ldots,
\bar x^{R_p}$.
\begin{defis}
\begin{enumerate}
\item We shall say that the family $M$ of manifolds is {\bf non-resonnant} whenever, for all $1\leq i\leq m$,  $1\leq k\leq n$ and for all $Q\in \Bbb N_2^n$, there exists a $1\leq j\leq m$ such that $\bar\mu_{i,j}^Q\neq \mu_{i,j,k}^{-1}$.
\item We shall say that the family $M$ of manifolds {\bf non-resonnant on ${\cal I}$} whenever for all monomial $z^Q$ not belonging to ${\cal I}$ and for all couple $(i,k)$, there exists $j$ such that $\bar\mu_{i,j}^Q\neq \mu_{i,j,k}^{-1}$.
%\item The $\rho_i$'s are simultaneously {\bf normalizable on $V({\cal I})$} whenever $R_i(\bar z) - D_{i,j}R_i(\bar %D_{i,j}\bar z)\in {\cal I}$ for all $1\leq i,j\leq l$. In this case, we have
%$$
%\rho_{i,k}(z)=B_i\bar z + \sum_{\substack{Q\in \Bbb N^n_2, z^Q\not\in {\cal I}\\ \bar\mu_{i,j}^Q=\mu_{i,j,k}^{-1}}} %\rho_{i,k,Q}\bar z^Q.
%$$
\end{enumerate}
\end{defis}
\begin{theos}\label{theo-main}
Assume that the group $G$ is abelian. Let ${\cal I}$ be
a monomial ideal (resp. properly imbedded) left invariant by the family $D:=\{D_{i,j}\}$ and the $B_i$'s . Assume that $D$ is diophantine (resp. on ${\cal I}$) and that $M$ is non-resonnant on ${\cal I}$. Assume $G$ is formally linearizable on ${\cal I}$.
Then, the family $F$ is holomorphically linearizable on ${\cal I}$. Moreover, in these coordinates, the $\rho_i$'s are anti-linearized on $\bar{\cal I}$.
\end{theos}
\begin{proof}
By theorem \ref{theo-invariant}, the family $F$ is holomorphically linearized on ${\cal I}$. Let us show that, in these coordinates, the $\rho_i$'s are anti-linearized on $\bar{\cal I}$.

Let us prove by induction on $|Q|\geq 2$ that $\{\rho_{i,k}\}_Q=0$ whenever
$z^Q$ doesn't belong to ${\cal I}$ and
$\bar\mu_{i,j}^Q\neq\mu_{i,j,k}^{-1}$. We recall that $\{\rho_{i,k}\}_Q$
denotes the coefficient of $\bar z^Q$ in the Taylor expansion of
$\rho_{i,k}$. Assume it is case up to order $k$. Let $Q\in \Bbb N_2^n$ with $|Q|=k+1$. Let us compute $\{\rho_{i,k}\}_Q$. Using equation $(\ref{conj})$, we obtain
\begin{eqnarray*}
R_{i}(\bar z) - D_{i,j}R_{i}(\bar D_{i,j}\bar z)& = & D_{i,j}B_i\bar f_{i,j}(\bar z)+f_{i,j}(B_j\bar z)\\
+D_{i,j}\left(R_{i}(\bar D_{i,j}\bar z) -R_{i}(\bar F_{i,j}\bar z)\right)& & +\left(f_{i,j}(\rho_j)-f_{i,j}(B_j\bar z)\right).
\end{eqnarray*}
Moreover, $F$ is linearized on $V({\cal I})$. Hence, both $\{D_{i,j}B_i\bar f_{i,j}(\bar z)+f_{i,j}(B_j\bar z)\}_Q$ and $\{R_{i}(\bar D_{i,j}\bar z) -R_{i}(\bar F_{i,j}\bar z)\}_Q$ vanish when $z^Q$ doesn't belong to ${\cal I}$. Hence, if $z^Q\not\in {\cal I}$, then we have
$$
(1-\mu_{i,j,k}\bar\mu^Q_{i,j})R_{Q,i,k}=\{\left(f_{i,j,k}(\rho_j)-f_{i,j,k}(B_j\bar z)\right\}_Q.
$$
But by induction, we have 
$$
\{\left(f_{i,j,k}(\rho_j)-f_{i,j,k}(B_j\bar z)\right\}_Q = \{Df_{i,j,k}(B_j\bar z)R_j+Df_{i,j,k}^2(B_j\bar z)R_j^2+\cdots\}_Q =0.
$$
Therefore, since $(1-\mu_{i,j,k}\bar\mu^Q_{i,j})\neq 0$, then we have $R_{Q,i,k}=0$. That is, 
$$
\rho_i(z) = B_i\bar z \mod \bar{\cal I}.
$$
\end{proof}
\begin{coros}
Under the assumptions of theorem \ref{theo-main}, there exists a complex analytic subvariety ${\cal S}$ passing throught the origin and intersecting each totally real submanifold $M_i$. In good holomorphic coordinate system, ${\cal S}$ is a finite intersection of a finite union of complex hyperplane defined by complex coordinate subspaces :
$$
{\cal S}=\cap_i\cup_j\{z_{i_j}=0\}.
$$
The intersection $M_k\cap {\cal S}$ is then given by 
$$
M_k\cap {\cal S} = \left\{z\in \cap_i\cup_j\{z_{i_j}=0\}\,|\; B_k\bar z =z\right\}.
$$
\end{coros}
\begin{proof}
The complex analytic subvariety ${\cal S}$ is nothing but $V({\cal I})$. The trace of it on $M_i$ is the fixed points set of $\rho_i$ belonging to $V({\cal I})$. It is non void since it contains the origin. According to the previous theorem, the $\rho_i$'s are holomorphically and simultaneously linearizable on $V({\cal I})$. By assumptions, ${\cal I}$ is a monomial ideal so $V({\cal I})$ is a finite intersection of a finite union of hyperplane defined by coordinate subspaces :
$$
{\cal S}=\cap_i\cup_j\{z_{i_j}=0\}.
$$
 
\end{proof}
\begin{coros}
Assume that the family $M$ is non-resonnant, $G$ is formally linearizable and $D$ is diophantine. Then, in a good holomorphic coordinates system, $M$ is composed of linear totally real subspaces
$$
\bigcup_i\left\{z\in\Bbb C^n\,|\; B_i\bar z = z\right\}.
$$
\end{coros}
\begin{rems}
If the family $M$ is non-resonnant and if for all $(i,k)$, one of the eigenvalues $\mu_{i,j,k}$'s belong to the unit circle, then $G$ is formally linearizable. In fact, for any $Q\in \Bbb N_2^n$, any $1\leq i\leq m$, any $1\leq k\leq n$, there exists $1\leq j\leq m$ such that 
$$
\bar\mu_{i,j}^Q \neq \mu_{i,j,k}^{-1}=\bar \mu_{i,j,k}.
$$
This means precisely that $D$ is non-resonnant in the classical sense. There is no obstruction to formal linearization.
\end{rems}
\begin{coros}
Let ${\cal I}$ be the ideal generated by the monomials $x^{R_1},\ldots,
x^{R_p}$ generating the ring $\widehat{\cal O}_n^D$ of formal invariants of
$D$. We assume that the non-linear centralizer of $D$ is generated by the same
monomials. If $D$ is diophantine on ${\cal I}$ then, in a good holomorphic coordinate system, we have 
$$
V({\cal I})=\{z\in (\Bbb C^n,0)\;|\; z^{R_1}=\cdots=z^{R_p}=0 \},
$$
and 
$$
\rho_{i|V({\cal I})}(z)= B_i\bar z.
$$
\end{coros}
\begin{coros}
Let us consider two totally real $n$-manifolds of $(\Bbb C^n,0)$ not
intersecting transversally at the origin. Assume that the $l$ first
eignevalues of $DF(0)$ are one. Let $\mu^{R_1}=1,\ldots, \mu^{R_p}=1$ be the
other (i.e. $R_i\in \Bbb N^n$ and $|R_i|>1$) generators of resonnant
relations. 
Let 
$$
V({\cal I})=\{z\in (\Bbb C^n,0)\;|\; z_1=\cdots =z_l=z^{R_1}=\cdots=z^{R_p}=0 \}.
$$
If $DF(0)$ is diophantine on $V({\cal I})$, then in good holomorphic
coordinate system,
$$
M_i\cap V({\cal I}) = \left \{z\in V({\cal I})\,|\; B_i\bar z = z\right\},\quad i=1,2.
$$
\end{coros}

%\subsection{When the associated group is not abelian}

%%% Local Variables: 
%%% mode: latex
%%% TeX-master: "webster"
%%% End: 

\section{Real analytic manifolds with CR singularities}

Les us consider a $(n+p-1)$-real analytic submanifold $M$ of $\Bbb C^n$ of the form
\begin{equation}\label{variete-orig}
\begin{cases} z_1 = x_1+iy_1\\ \vdots \\ z_p = x_p+iy_p\\ 
y_{p+1} = F_{p+1}(z',\bar z', x'')\\ \vdots\\ y_{n-1} = F_{n-1}(z',\bar z', x'')\\
z_{n} = G(z',\bar z', x'')\end{cases}
\end{equation}
where we have set $z'=(z_1,\ldots,z_p)$, $z''=(z_{p+1},\ldots, z_{n-1})$.
The real analytic real (resp. complex) valued functions $F_i$ (resp. $G$) are assumed to vanish at the origin as well as its derivative. The tangent space at the origin contains the complex subspace defined by $z'$.

First of all, we shall show, under some assumptions, that there is a good
holomorphic coordinate systems in which the $F_i$'s are of order greater than
or equal to $3$ and the $2$-jet of the $G$'s depends only on $z'$ and its
conjugate. When $p=1$, this was done by E. Bishop ($n=2$)\cite{bishop},
and also by J. Moser and S. Webster ($n\geq 2$)\cite{moser-webster}.

\subsection{Preparation}

Let us consider the $2$-jet $G^2$ of $G$. It can be written as a sum of a quadratic polynomial 
$$
Q(z',\bar z'):= \sum d_{i,l} z'_i\bar z'_l + \sum_{1\leq i,l\leq p} e_{i,l} z'_i z'_l +
\sum_{1\leq i,l\leq p} f_{i,l} \bar z'_i\bar z'_l
$$ 
and 
$$
\Sigma:=\sum a_{\alpha,\beta}x''_{\alpha}x''_{\beta} + \sum b_{\alpha,i}x''_{\alpha}z'_{i} +\sum c_{\alpha,i}x''_{\alpha}\bar z'_{i}.
$$

\begin{defis}
Let 
$$
 F:=\left(\begin{matrix}f_{1,1}& \ldots & f_{p,1} \\ \vdots & & \vdots \\ f_{1,p}& \ldots & f_{p,p} \end{matrix}\right)
$$
be the matrix associated to the so normalized $Q$. The symmetric matrix $S=\frac{1}{2}(F+ ^tF)$ will be called a {\bf Bishop matrix}.
\end{defis}

First of all, by a linear change of the coordinates $z'$, we can diagonalize
the symmetric matrix $S$.

Our {\bf first non-degeneracy condition} is that the sesquilinear part $\tilde
H:=\sum d_{i,l}
z'_i\bar z'_l$ of $Q$ is not reduced to $0$.
Let 
$$
\|\tilde H\|:= \sup_{z\in\Bbb C^p\,|\,\|z\|=1} \frac{|\sum d_{j,i,l} z'_i\bar z'_l|}{\|z'\|^2}\neq 0.
$$
to its norm. Let us set $Z_n=z_n/\|\tilde H\|$.
Then the sequilinear part $H$ of the new $Q$ is of norm $1$. By setting 
$$
Z_n:= z_n + \sum (f_{j,i,l}-e_{j,i,l}) z'_i z'_l,
$$
we transform $Q$ in the following form~:
\begin{equation}\label{Q}
Q(z',\bar z') = H(z,'\bar z')+ \sum_{i=1}^p  \gamma_{i}( (z'_i)^2+ (\bar z'_i)^2).
\end{equation}
\begin{defis}
The eigenvalues $\gamma_{1},\ldots, \gamma_{n}$ of the so normalized Bishop matrix $S$ will be called the {\bf generalized Bishop invariants}.
\end{defis}

Let us show that, by an holomorphic change of coordinates, we can get rid of $\Sigma$. First of all, let us get rid of the third member of $\Sigma$ by a change of the form
$$
\zeta'_i \mapsto z'_i:=\zeta'_i + \sum_{\gamma=p+1}^{n-1} A_{i,\gamma} z''_{\gamma},\quad i=1,\ldots, p.
$$
We have
\begin{eqnarray*}
G(z',\bar z', x'') & =  & \sum a_{\alpha,\beta}x''_{\alpha}x''_{\beta} + \sum b_{\alpha,i}x''_{\alpha}\left(\zeta'_i + \sum_{\gamma=p+1}^{n-1} A_{i,\gamma}z''_{\gamma}\right)\\
& &+\sum c_{\alpha,i}x''_{\alpha}\left(\zeta'_i + \sum_{\gamma=p+1}^{n-1} \bar A_{i,\gamma}\bar z''_{\gamma}\right) \\
& &  +\sum d_{i,l}\left(\zeta'_i + \sum_{\gamma=p+1}^{n-1} A_{i,\gamma}z''_{\gamma}\right)\left(\bar \zeta'_l + \sum_{\gamma=p+1}^{n-1} \bar A_{l,\gamma}\bar z''_{\gamma}\right)\\
&& + \sum f_{i,l} \left(\zeta'_i + \sum_{\gamma=p+1}^{n-1} A_{i,\gamma}z''_{\gamma}\right) \left(\zeta'_l + \sum_{\gamma=p+1}^{n-1} A_{l,\gamma}z''_{\gamma}\right) \\
&&+ \sum f_{i,l} \left(\bar \zeta'_i + \sum_{\gamma=p+1}^{n-1} \bar
  A_{i,\gamma}\bar z''_{\gamma}\right)\left(\bar \zeta'_l +
  \sum_{\gamma'=p+1}^{n-1} \bar A_{l,\gamma'}\bar z''_{\gamma'}\right)\\
&& + \text{higher order terms.}
\end{eqnarray*}
The coefficient of $x_{\alpha}''\bar \zeta_i$ is
$$
c_{\alpha,i}+\sum_s d_{s,i} A_{s,\alpha}+\sum_s (f_{i,s}+f_{s,i})\bar A_{s,\alpha}
$$
We assume that we can solve the set of equations
\begin{eqnarray*}
-c_{\alpha,i}& = & \sum_s d_{s,i} A_{s,\alpha}+\sum_s (f_{i,s}+f_{s,i})\bar A_{s,\alpha}\\
-\bar c_{\alpha,i} & = & \sum_s \bar d_{s,i} \bar A_{s,\alpha}+\sum_s (\bar f_{i,s}+\bar f_{s,i}) A_{s,\alpha}.
\end{eqnarray*}
For each $p+1\leq \alpha\leq n-1$, let us set 
$$
A_{\alpha}:=\left(\begin{matrix}A_{1,\alpha}\\ \vdots \\ A_{p,\alpha}\end{matrix}\right),\; C_{\alpha} := \left(\begin{matrix}c_{\alpha,1}\\ \vdots \\ c_{\alpha,p}\end{matrix}\right),\;D:=\left(\begin{matrix}d_{1,1}& \ldots & d_{p,1} \\ \vdots & & \vdots \\ d_{1,p}& \ldots & d_{p,p} \end{matrix}\right).
$$
The previous equations can be written as
\begin{equation}\label{condition1}
\begin{cases}
-C_{\alpha}  =  DA_{\alpha}+(F+ ^tF)\bar A_{\alpha}\\
-\bar C_{\alpha}  =  \bar D\bar A_{\alpha}+(\bar F+ ^t\bar F)A_{\alpha}
\end{cases}
\end{equation}
These systems can be solved whenever
%\begin{equation}\label{cond1}
%\det \left(D\bar D-4\bar S S\right)\neq 0,
%\end{equation}
\begin{equation}\label{cond1}
\det \left(4^{-1}S^{-1}D\bar S^{-1}\bar D-I_p\right)\neq 0,
\end{equation}
where $S:=1/2(F+ ^tF)$ is the Bishop matrix.
\begin{rems}
When $p=1$, this condition reads $4\gamma^2\neq 1$ where $\gamma $ is the Bishop invariant \cite{bishop}.
\end{rems}

Now let us remove the two first terms of the $\Sigma$'s (once we have done the previous change of coordinates, the coefficients appearing in $\Sigma$ may have changed). Let us set
$$
Z_n:= z_n-\sum a_{\alpha,\beta}z''_{\alpha}z''_{\beta} - \sum b_{\alpha,i}z''_{\alpha}z'_{i}.
$$
Hence, we have
\begin{eqnarray*}
Z_n & = & G(z',\bar z', x'')-\sum a_{\alpha,\beta}z''_{\alpha}z''_{\beta} - \sum b_{\alpha,i}z''_{\alpha}z'_{i} \\
& = & Q(z',\bar z') +\Sigma -\sum a_{\alpha,\beta}z''_{\alpha}z''_{\beta} -
\sum b_{\alpha,i}z''_{\alpha}z'_{i}\\
& & + \text{higher order terms}.
\end{eqnarray*}
But we have 
$$
z''_{\alpha}=x''_{\alpha}+iF_{\alpha}(z',\bar z', x'')
$$
where $F_{\alpha}$ is of order greater than or equal to two. Hence, the the 2-jet of $Z_n$ is precisely $Q(z',\bar z', x'')$.

Let us consider the $F_{\alpha}$'s. As above, the $2$-jet of $F_{\alpha}$ can be written as the sum of 
$$
Q_{\alpha}(z',\bar z'):= \sum d_{\alpha,i,l} z'_i\bar z'_l + \sum e_{\alpha,i,l} z'_i z'_l + \sum \bar e_{\alpha,i,l} \bar z'_i\bar z'_l
$$ 
and 
$$
\Sigma_{\alpha}:=\sum a_{\alpha,\gamma,\beta}x''_{\gamma}x''_{\beta} + 2\text{Re}\sum b_{\alpha,\gamma,i}x''_{\gamma}z'_{i},
$$
where $a_{\alpha,\gamma,\beta}$ is a real number, $d_{\alpha,i,l}= \bar d_{\alpha,l,i}$. 
Let us show that, under some other assumption, we can get rid of the $z_i'\bar z_j'$'s terms in the $Q_{\alpha}'s$, $\alpha=p+1,\ldots, n-1$.
In fact, let us set
$$
z_{\alpha}\mapsto Z_{\alpha} := z_{\alpha} + i b_{\alpha}z_n''',\quad \alpha=p+1,\ldots, n-1.
$$
In the new coordinates, the $2$-jet of $Z_{\alpha}$ is 
$$
iQ_{\alpha}(z',\bar z')+ib_{\alpha}Q(z',\bar z')+i\Sigma_{\alpha}.
$$
Thus, in order that the coefficients of the $z_i'\bar z_j'$'s vanish, one must
satisfies the {\bf second non-degeneracy condition}~:\begin{equation}\label{condition2} \text{{\it each $Q_{\alpha}$ is proportional to $Q$.}}\end{equation}

Let us show that we can get rid of the quadratic part of $F_{\alpha}$ by the using following change of coordinates :
$$
z_{\alpha}\mapsto Z_{\alpha}:= z_{\alpha}-2i\sum b_{\alpha,\gamma,i}z''_{\gamma}z'_{i}-i\sum a_{\alpha,\gamma,\beta}z''_{\gamma}z''_{\beta}-2i\sum e_{\alpha,i,l} z'_i z'_l. 
$$
In fact, we have
$$
\text{Im}(Z_{\alpha})  =  \text{Im}\left( z_{\alpha}-2i\sum b_{\alpha,\gamma,i}z''_{\gamma}z'_{i} -i\sum a_{\alpha,\gamma,\beta}z''_{\gamma}z''_{\beta}-2i\sum e_{\alpha,i,l} z'_i z'_l\right) 
$$
The $2$-jet of the right hand side is 
$$
Q_{\alpha}(z',\bar z', X'')+ \Sigma_{\alpha}(z',\bar z', X'')+ \text{Im}\left( -2i\sum b_{\alpha,\gamma,i}X''_{\gamma}z'_{i} -i\sum a_{\alpha,\gamma,\beta}X''_{\gamma}X''_{\beta}-2i\sum e_{\alpha,i,l} z'_i z'_l\right).
$$
It is equal to 
$$
\sum d_{\alpha,i,l} z'_i\bar z'_l
$$
since
$$
\text{Im}\left(2i\sum e_{\alpha,i,l} z'_i z'_l\right)= 2\text{Re}\sum e_{\alpha,i,l} z'_i z'_l.
$$
We summarize these results in the following lemma.
\begin{lemms}\label{bishop-gen}
Under assumptions $(\ref{cond1})$ and $(\ref{condition2})$, the manifold $(\ref{variete-orig})$ can be transformed, by an holomorphic change of variables, to 
\begin{equation}\label{variete-arriv}
\begin{cases} z_1 = x_1+iy_1\\ \vdots \\ z_p = x_p+iy_p\\ 
y_{p+1} =  f_{p+1}(z',\bar z', x'')\\ \vdots\\ y_{n-1} =  f_{n-1}(z',\bar z', x'')\\
z_{n} = Q(z',\bar z')+ g(z',\bar z', x'')\end{cases}
\end{equation}
where the $f_i$'s and $g$ are germs of real analytic functions at the origin
and of order greater than or equal to $3$ there. The quadratic polynomial $Q$ is of the form 
$$
Q(z',\bar z') = \sum d_{i,l} z'_i \bar z'_l+ \sum  \gamma_{i}( (z'_i)^2
+(\bar z'_i)^2),
$$
the norm of the sesquilinear part of $Q$ being $1$.
\end{lemms}

In the next two sections, we shall adapt the constuction of Moser and Webster to our context.
\subsection{Complexification}

Let us complexify such a manifold $M$ by replacing $\bar z_i$ by $w_i$ in order
to obtain a complex analytic $(n+p-1)$-manifold ${\cal M}$ of $\Bbb C^{2n}$: 
\begin{equation}\label{variete-complex}
\begin{cases} 
2x_{\alpha} = z_{\alpha}+w_{\alpha},\\
z_{\alpha}-w_{\alpha} = 2iF_{\alpha}(z',w', x'')= 2i\bar F_{\alpha}(w',z', x''),\\
z_{n} = G(z',w', x''),\\
w_{n} = \bar G(w',z', x''),\\
\end{cases}
\end{equation}
where $\alpha$ ranges from $p+1$ to $n- 1$. In this situation, we have $G(z',w',x")=Q(z',w')+g (z',w',x")$ where $g$ as well as the $F_{\alpha}$'s are of order greater than or equal to $3$. As usual if $G(y)=\sum_{Q}g_{Q}y^Q$ is a formal power series in $\Bbb C^k$,
then $\bar G$ denotes the formal power series $\sum_{Q}\bar g_{Q}y^Q$.
These equations imply that 
$$
z_{\alpha}= x_{\alpha}+iF_{\alpha}(z',w', x''); \quad w_{\alpha}= x_{\alpha}-iF_{\alpha}(z',w', x'').
$$
The variables $(z',w',x'')$ may be used as complex coordinates.
Let us define the anti-holomorphic involution $\rho$ of $\Bbb C^{2n}$ by
$$
\rho(z,w)=(\bar w,\bar z).
$$
A complex anaytic manifold ${\cal M}$ of $\Bbb C^{2n}$ comes from a manifold
$M$ whenever it is preserved by $\rho$. In this case, $M={\cal M}\cap
\text{Fix}(\rho)$. The restriction to ${\cal M}$ of this map is the anti-holomorphic involution $\rho(z',w',x'')=(\bar w',\bar z',\bar x'')$.

The two projections $\pi_1(z,w)=z$ and $\pi_2(z,w)=w$, when restricted to ${\cal M}$ have the form
\begin{eqnarray*}
\pi_1(z,w)& = & (z', x_{\alpha}+iF_{\alpha}(z',w', x''), G(z',w', x''))\\
\pi_2(z,w)& = & (w', x_{\alpha}-iF_{\alpha}(z',w', x''), \bar G(w', z', x'')).
\end{eqnarray*}

Let us define the holomorphic involution $\tau_1(z,w)=(\tilde z,\tilde w)$
(resp.$\tau_2(z,w)=(\tilde z,\tilde w)$)  on ${\cal M}$ by $w=\tilde w$
(resp. $z=\tilde z$). This leads to the following equations :
\begin{eqnarray*}
\tilde w' & = & w'\\
\tilde x_{\alpha}+iF_{\alpha}(\tilde z',w', \tilde x'') & = & x_{\alpha}+iF_{\alpha}(z',w', x''),\quad \alpha=p+1,\ldots, n-1\\
\bar G(w',\tilde z', \tilde x'') & = & \bar G(w', z', x'').
\end{eqnarray*}

By the implicit function theorem, there exist an holomorphic function
$\Gamma_{\alpha}$ such that 
$$
\tilde x_{\alpha}= \Gamma_{\alpha}(z',\tilde z',w', \tilde x'')= x_{\alpha}\;
\text{mod}\; {\cal M}^{2}\quad \alpha=p+1,\ldots, n-1.
$$
Here, ${\cal M}$ denotes the maximal ideal of germ of holomorphic functions in
$\Bbb C^{n+2p-1}$ at $0$.
Hence, we have
\begin{equation}\label{equation-involution}
\bar G(w',\tilde z', \Gamma''(z',\tilde z',w', \tilde x''))  =  \bar G(w', z', x'').
\end{equation}

{\bf Assume that the Bishop matrix is invertible.}

%$$
%\tilde z' = \hat H(z', w', x'') = Mz'\; \text{mod}\;\widehat{\cal M}^2 
%$$
%where $M$ is a nonzero matrix and $\widehat{\cal M}$ is the maximal ideal of
%formal power series of $\Bbb C^{2p+q}$. Let $\hat \tau_1$ be the associated formal
%solution map. In order to compute $M$, let us solve

Let us solve, in $\tilde z'$, the following equation~:
$$
Q(z',w')=Q(\tilde z',w'),
$$
where $Q$ is the quadratic part of $G$. Therefore, we have
$$
Q(\tilde z',w')-Q(z',w') = D_zQ(z',w')(\tilde z'- z')+\frac{1}{2}D^2_zQ(z',w')(\tilde z' - z')^2.
$$
Since $\frac{\partial^2 Q}{\partial z_j\partial z_k}=0$ if $j\neq k$, we
obtain that, for any $i_0$, 
$$
\frac{\partial Q}{\partial z_{i_0}}(z',w') + \frac{1}{2}\frac{\partial^2 Q}{\partial z_{i_0}^2}(z',w')(\tilde z_{i_0} - z_{i_0})=0.
$$

But, for any $1\leq i_0\leq p$, we $1/2\frac{\partial^2 Q}{\partial z_{i_0}^2}(w',z') =
\gamma_{i_0}\neq 0$. Hence, we have
$$
\tilde z_{i_0} = -  z_{i_0} -\frac{1}{\gamma_{i_0}}\sum_{l=1}^pd_{i_0,l}w_l.
$$ 

Hence, by the implicit function theorem, equation $(\ref{equation-involution})$ admit an holomorphic solution $\tau_1(z',w', x'')$
which linear part at the origin is precisely $T_1z'$~:
$$
\tilde z' = H(z', w', x'')= T_1z'\; \text{mod}\;{\cal M}^2,\quad \tilde w'= w',\quad \tilde x''= \Gamma (z', w', x'').
$$
Here ${\cal M}$ denotes the maximal ideal of germs of holomorphic functions of
$(\Bbb C^{n+p-1},0)$ and the linear part $T_1$ is
$$
T_1:=D\tau_1(0)=\left(\begin{matrix}-Id_p & - S^{-1}D & 0 \\ 0 & Id_p & 0 \\ 0 &
    0 & Id_q \end{matrix}\right),
$$
where $D$ is the matrix $(d_{i,j})_{1\leq i,j\leq p}$ and $I_q$ stands for the
identity matrix of dimension $q:=n-p-1$. The map $\pi_2$ is a two-fold covering with covering
transformation $\tau_1$.

Since $\tau_2= \rho\circ\tau_1\circ\rho$, $\tau_2$ restricted to ${\cal M}$ is defined by
$$
\tilde z' = z',\quad \tilde w'= \bar H(w', z', x''),\quad \tilde x''= \bar
\Gamma (w', z', x'').
$$
The linear part at the origin of $\tau_2$ is 
$$
T_2:=\left(\begin{matrix}Id_p & 0 & 0 \\ - \bar S^{-1}\bar D & -Id_p & 0 \\ 0 &
    0 & Id_q \end{matrix}\right).
$$
The map $\pi_1$ is a two-fold covering with covering
transformation $\tau_2$.

Let us consider the germ of holomorphic diffeomorphism $g$ of $(\Bbb
C^{n+p-1},0)$ defined to be $g:=\tau_1\circ \tau_2$. It fixes the origin. The linear part of the
diffeomorphism $g$ at the origin is 
$$
\Phi:=Dg(0)=T_1T_2=\left(\begin{matrix}-Id_p +S^{-1}D\overline{S^{-1}D}  & S^{-1}D & 0 \\ -\overline{S^{-1}D} & -I_p & 0\\0 & 0 &
    I_q\end{matrix}\right).
$$

\subsection{Quadrics and linear involutions}

In this section will shall study the relation between the linear part of the
involutions, the linear anti-holomorphic involution and the quadric.

Under our assumptions, the set of fixed points of $T_1$ and $T_2$ is a
$q$-dimensional vector space ($q=n-p-1$) and these are the only common eigenvectors. In
fact, let us try to solve $T_1v=av$ and $T_2v=bv$. Let us write
$v=(v_1,v_2,v_3)$ with $v_1,v_2$ belong to $\Bbb C^p$ while $v_3$ belongs to
$\Bbb C^q$. This leads to
$$
(*)\left\{\begin{array}{l}-v_1-S^{-1}Dv_2=av_1\\ v_2 = a v_2\\v_3 = a
    v_3\end{array}\right.\quad (**)\left\{\begin{array}{l}v_1=bv_1\\ -\bar
    S^{-1}\bar Dv_1 -v_2= b v_2\\v_3 = b v_3\end{array}\right.
$$
Assume $v_3=0$. If $v_2=0$ then $a=-1$ (otherwise $v_1=0$ too). The second
equation of $(**)$ gives $\bar S^{-1}\bar Dv_1=0$; that is $v_1=0$. This is
not possible. Thus $v_2\neq 0$ and $a=1$. The first equation of $(*)$ gives
$-S^{-1} Dv_2=2v_1$. Moreover, we have $b=1$ since otherwise we would have
$v_1=0$ and $v_2=0$ by the second equation of $(**)$. Therefore, using the
same equation we obtain 
$$
\bar S^{-1}\bar D S^{-1} Dv_2= 4v_2.
$$
According to condition $(\ref{cond1})$, we have $v_2=0$ and then $v_1=0$. Now, if $v_3\neq 0$
then $a=b=1$, then we can apply the previous resonning to obtain $v_1=v_2=0$.

Let $V_i$ be the $(-1)$-eigenspace of $T_i$, $i=1,2$. Let $E$ be their common
$(+1)$-eigenspace. {\bf We assume that $V_1,V_2$ and $E$ span $\Bbb C^{n+p-1}$}. Let
$F=V_1+V_2$. We have $\Bbb C^{n+p-1}=F\oplus E$. Let us show that it is invariant under both $T_1$ and $T_2$. In fact,
let $v_2\in V_2$. A priori, we have $T_1v_2 = \tilde v_1 + \tilde v_2 + e$ where
$\tilde v_1$ (resp. $\tilde v_2$, $e$) belongs to $V_1$ (resp. $V_2$, $E$). Since
$T_1^2=Id$ and $T_1v_1=-v_1$, we have $v_2 = -\tilde v_1 + T_1\tilde v_2 +
e$. But, $T_1\tilde v_2= \tilde v'_1 + \tilde v'_2 + e'$. Hence, we have $v_2
= -\tilde v_1 +\tilde v'_1 + \tilde v'_2 + e'+ e$. Therefore, $-\tilde v_1
+\tilde v'_1=0$, $\tilde v'_2-v_2=0$ and $e'+ e=0$. It comes
$$
T_1(v_2-\tilde v_2-e)= - (v_2-\tilde v_2-e).
$$
This means that $v_2-\tilde v_2-e$ belongs to $V_1$ and leads to $e=0$; that
is $F$ is left invariant by $T_1$. The same argument applies to $T_2$.
Let $T'_1$ (resp. $T'_2$, $\Phi'$) be the restriction $T_1$ (resp. $T_2$,
$\Phi$) to $F$.

Let $\mu$ be an eigenvalue of $\Phi'$ corresponding to an eigenvector
$f$. Then, $T'_2f= \mu T'_1 f = \mu\Phi' T'_2 f$. Hence, $\mu^{-1}$ is another
eigenvalue of $\Phi'$ associated to $T'_2f$. The vectors $f$ and $T'_2f$ are independant
since otherwise we would have $cf=T'_2f= \mu T'_1 f$ and $T'_1$ and $T'_2$
would have a common eigenvector. If $\mu=\mu^{-1}$, then the restriction of
$\Phi'$ to the span $V_f$ of $f$ and $T'_2f$ is $\pm Id$. This implies a common
eigenvector to $ T'_2$ and $ T'_1$. As in \cite{moser-webster}[p.269], we can
choose $\{f,T'_2f \}$ as a basis of $V_f$. In this basis, the matrix of $T'_1$, $T'_2$ and $\Phi'$
are
$$
\Phi'_{|V_f}=\left(\begin{matrix}\mu & 0 \\ 0 &
    \mu^{-1}\end{matrix}\right)\quad T'_{i|V_f}=\left(\begin{matrix}0& \lambda_i \\ \lambda_i^{-1} & 0\end{matrix}\right),
$$
where $\lambda_1=\mu$ and $\lambda_2=1$.

Let us show that $E$ is invariant under the linear anti-holomorphic involution $\rho$. Let $e$ belongs to $E$. By definition, it is left invariant by both $T_1$ and $T_2$. Since, $T_1\rho =\rho T_2$, we obtain $T_1\rho(e)=\rho(T_2e)=\rho(e)$ and $T_2\rho(e)=\rho(e)$. Hence, $\rho(E)=E$. 

Let us show that $F$ is invariant under the linear anti-holomorphic involution $\rho$. 
Let $N$ be the totally real fixed point set of $\rho$ on $E$. This means that $E=N+iN$ and that we may choose coordinates $\zeta$ on $E$ so that $\rho :\zeta \mapsto \bar \zeta$. Let us show that $\rho(F)$ is also invariant by both $T_1$ and $T_2$. In fact, we have $T_1\rho(F) =\rho T_2(F)=\rho(F)$ and $T_2\rho(F) =\rho(F)$ as well. Hence, $E\oplus\rho(F)=\Bbb C^{n+p-1}$ is a decomposition preserved by both $T_1$ and $T_2$. The space $\rho(F)$ has to contain the $(-1)$-eigenspace of both $T_1$ and $T_2$, that is $F$.

Let $f$ be a $\mu$-eigenvector of $\Phi$. As above, $T_2f$ is a
$\mu^{-1}$-eigenvector of $\Phi$. Since $\Phi \rho\Phi= T_1\rho T_2=\rho$, we have
$$
\rho(f) = \bar\mu \Phi(\rho(f)).
$$
Hence, $\rho(f)$ is a ${\bar\mu}^{-1}$-eigenvector of $\Phi$. Moreover,
$T_2\rho(f)=\rho(T_1f)=\bar\mu^{-1} \rho(T_2f)$ is a ${\bar\mu}$-eigenvector
of $\Phi$. This means that $\rho(V_f)=V_{\rho(f)}$.

Let us assume that the set ${\bar\mu}^{-1}$ is different than
${\mu}^{-1}$ and $\mu$. Then, $\rho(T_2\rho(f))=T_1f$ is a
${\mu}^{-1}$-eigenvector of $\Phi$. So, $T_2\rho(T_2\rho(f))=T_2T_1f=\Phi^{-1}(f)=\mu^{-1}f$ is a
${\mu}$-eigenvector of $\Phi$. In this case, the matrices of $\rho$ and $\Phi$
restricted to $V_f\oplus V_{\rho(f)}$ are
$$
\rho_{|V_f\oplus V_{\rho(f)}}=\left(\begin{matrix}0 & 0 & 1 & 0 \\ 0 & 0 & 0
    &\mu^{-1}\\ 1 & 0 & 0 & 0 \\ 0 & \bar \mu & 0 & 0 \end{matrix}\right)\quad
\Phi_{
|V_f\oplus V_{\rho(f)}}=\left(\begin{matrix}\mu & 0 & 0 & 0\\ 0 & \mu^{-1} & 0
  & 0 \\0 & 0 & \bar \mu^{-1} & 0\\0 & 0 & 0 & \bar \mu\end{matrix}\right).
$$
%Therefore, we have ${\bar\mu}^{-1}={\mu}^{-1}$ or${\bar\mu}^{-1}={\mu}$. 
We recall that $V_f$ denotes the span of $f$ and $T_2f$.
We can do a similar analysis as in Moser-Webster article \cite{moser-webster}[p.269]~: let $\mu_i$ be
an eigenvalue of $\Phi$ of multiplicity $n(i)$ such that $\mu_i=\bar
\mu_i^{-1}$ or $\mu_i=\bar\mu_i$. Let $\{f_1,\ldots, f_{n(i)}\}$
  be an associated basis. Let $E_i$ be the span of the basis $\{f_1,\ldots,
  f_{n(i)},T_2f_1,\ldots, T_2f_{n(i)}\}$. We have
\begin{itemize}
\item if $\mu_i=\bar\mu_i^{-1}$, then $\rho(f_k)=\sum a_{j,k}f_k$ so that 
$$
\rho(T_2f_k)=T_1\rho(f_k)=T_1(\sum a_{j,k}f_k)=\sum
a_{j,k}T_1f_k=\bar\mu_i\sum a_{j,k}T_2(f_k).
$$
\item if $\mu_i=\bar\mu_i$, then $\rho(f_k)=\sum a_{j,k}T_2f_k$ so that 
$$
\rho(T_2f_k)=T_1\rho(f_k)=T_1(\sum a_{j,k}T_2f_k)=\mu_i\sum a_{j,k}f_k.
$$
\end{itemize}
We have used the property that, if $\Phi(f)=\mu f$ then $T_2f=\mu T_1f$.
Hence we have proved the following 
\begin{lemms}\label{rho}
Let $\Phi$, $T_1$, $T_2$ and $\rho$ as above. Then there exists a decomposition
$$
\Bbb C^{n+p-1} = E_1\oplus\cdots\oplus E_r\oplus E_{r+1}\oplus\cdots\oplus E_{s}\oplus G
$$
left invariant by $\Phi$, $T_1$, $T_2$ and $\rho$ such that
\begin{itemize}
\item The $E_i$'s are complex vectorspaces of
  dimension $2n(i)$ if $1\leq i\leq r$  and $4n(i)$ otherwise. $G$ is a
  $n-p-1$-dimensional vectorspace.
\item The restrictions of $\Phi$, $T_1$, $T_2$ to $G$ is the identity
\item If $1\leq i\leq r$, then $\mu_i\in\{\bar\mu_i,\bar\mu_i^{-1}\}$ and there exists coordinates
  $(\zeta_{i},\eta_{i})$ of $E_i$ ($\zeta_{i}=(\zeta_{i,1},\ldots,\zeta_{i,n(i)})$) such that
$$
T_{k|E_i}=\left(\begin{matrix}0 &
    \lambda_{k,i}I_{n(i)}\\\lambda_{k,i}^{-1}I_{n(i)} & 0\end{matrix}\right),
\;k=1,2\quad \Phi_{|E_i}=\left(\begin{matrix}
    \mu_{i}I_{n(i)}& 0\\ 0 & \mu_{i}^{-1}I_{n(i)}0\end{matrix}\right),
$$
where $\lambda_{1,i}=\mu_i$ and $\lambda_{2,i}=1$.

\item If $0<j\leq s-r$, then $\mu_{r+2j-1}\not\in \{\bar\mu_{r+2j-1},\bar\mu_{r+2j-1}^{-1}\}$ and there
  exists coordinates $\theta_j:=(\zeta_{r+2j-1},\eta_{r+2j-1},\zeta_{r+2j},\eta_{r+2j})$ of
  $E_{r+j}$ such that 
$$
\Phi_{|E_{r+j}}=\left(\begin{matrix}
    \mu_{r+2j-1}I_{n(r+j)}& 0 & 0 & 0\\ 0 & \mu_{r+2j-1}^{-1}I_{n(r+j)} & 0 & 0 \\ 0 & 0 &
    \bar\mu_{r+2j-1}^{-1}I_{n(r+j)} & 0\\ 0 & 0 & 0 & \bar\mu_{r+2j-1}I_{n(r+j)}\end{matrix}\right).
$$
and 
$$
T_{k|E_{r+j}}=\left(\begin{matrix}0 &
    \lambda_{k,j+2r-1}I_{n(j+r)} & 0 & 0 \\\lambda_{k,j+2r-1}^{-1}I_{m(j+r)} & 0 & 0 & 0\\
    0 & 0 & 0 & \bar\lambda_{k,j+2r-1}^{-1}I_{m(j+r)}\\0 & 0 &
    \bar\lambda_{k,j+2r-1}I_{m(j+r)} & 0 \end{matrix}\right),
$$
where $k=1,2$, $\lambda_{1,r+2j-1}=\mu_{r+2j-1}$, $\mu_{r+2j}=\bar \mu_{r+2j-1}^{-1}$ and $\lambda_{2,r+2j-1}=1$.
\item If ${\bar\mu_i}={\mu_i}^{-1}$ (we say {\bf hyperbolic}),  
$\rho_{|E_i}(\zeta_i,\eta_i)= \left(\begin{matrix}A_i & 0\\ 0 & \bar\mu_i A_i \end{matrix}\right)\left(\begin{matrix}\bar\zeta_i\\ \bar\eta_i\end{matrix}\right)$
for some $n(i)$-square matrix $A_i$ such that $A_i\bar A_i=I_{n(i)}$.
\item If ${\bar\mu_i}={\mu_i}$ (we say {\bf elliptic}),  
$\rho_{|E_i}(\zeta_i,\eta_i)= \left(\begin{matrix} 0 &\mu_i A_i \\  A_i & 0 \end{matrix}\right)\left(\begin{matrix}\bar\zeta_i\\ \bar\eta_i\end{matrix}\right)$
for some $n(i)$-square matrix $A_i$ such that $\mu_iA_i\bar A_i=I_{n(i)}$.
\item If $\mu_{r+2j-1}\not\in\{\bar\mu_{r+2j-1},\bar\mu_{r+2j-1}^{-1}\}$, that is $0<j\leq s-r$, (we say {\bf complex}), 
$$
\rho_{|E_{r+2j-1}}(\theta_j)=\left(\begin{matrix}0 & 0 & I_{n(j+r)} & 0 \\ 0 & 0 & 0
    &\mu_{r+2j-1}^{-1}I_{n(j+r)}\\ I_{n(j+r)} & 0 & 0 & 0 \\ 0 & \bar \mu_{r+2j-1}I_{n(j+r)} & 0 & 0 \end{matrix}\right)\left(\begin{matrix}\bar\zeta_{r+2j-1}\\ \bar\eta_{r+2j-1}\\\bar\zeta_{r+2j}\\ \bar\eta_{r+2j}\end{matrix}\right)
$$
\item $\rho_{|G}(\upsilon)=\bar\upsilon$.
\end{itemize}
\end{lemms}
\begin{defis}
A coordinate or its associated eigenvalue will be called {\bf hyperbolic} (resp. {\bf
  elliptic}, {\bf complex}) if $\mu_i=\bar\mu_i^{-1}$,
  (resp. $\mu_i=\bar\mu_i$, $\mu_i\not\in\{\bar\mu_i,\bar\mu_i^{-1}\}$).
\end{defis}
\begin{rems}
When $p=1$, the complex case doesn't appear. In fact, in this case, the eigenspace
associated to an eigenvalue $\mu\neq 1$ is a $2$-dimensional vector
space. Hence, we must have $\mu\in\{\bar\mu,\bar\mu^{-1}\}$.
\end{rems}

\subsection{Submanifolds with  CR-singularities}

Let us consider a real analytic $(n+p-1)$-submanifolds $M$ of $\Bbb C^n$ passing through the origin of the form
\begin{equation}
\begin{cases} z_1 = x_1+iy_1\\ \vdots \\ z_p = x_p+iy_p\\ 
y_{p+1} =  F_{p+1}(z',\bar z', x'')\\ \vdots\\ y_{n-1} =  F_{n-1}(z',\bar z', x'')\\
z_{n} = G(z',\bar z', x''):=Q(z',\bar z')+ g(z',\bar z', x'')\end{cases}
\end{equation}
where the $F_{i}$'s and $g$ are real analytic in a neighbourhood of
the origin and of order greater than or equal to $3$ at $0$. We
assume that the quadratic polynomial $Q$ is normalized in the sense developped in the previous section. Hence, we have
$$
Q= \sum_{i,j}g_{i,j}z'_i\bar z'_j + \sum_{i=1}^p\gamma_{i}((z'_i)^2+(\bar z'_i)^2),
$$
where the sesquilinear part is of norm $1$, and the $\gamma_{i}$'s are the
generalized Bishop invariants. The CR-structure is singular at the
origin. For a $n$-variety, we have $Q=z_1\bar z_1 + \gamma(z_1^2+\bar z_1^2)$. Let us complexify the equations as in $(\ref{variete-complex})$:
\begin{equation*}
{\cal M}:\begin{cases} 
2x_{\alpha} = z_{\alpha}+w_{\alpha},\\
z_{\alpha}-w_{\alpha} = 2iF_{\alpha}(z',w', x'')= 2i\bar F_{\alpha}(w',z', x''),\\
z_{n} = G(z',w', x''),\\
w_{n} = \bar G(w',z', x''),\\
\end{cases}
\end{equation*}
where $\alpha$ ranges from $p+1$ to $n-1$.
Then there exists a pair $(\tau_{1},\tau_{2})$ of holomorphic involutions as
defined above. It defines a germ of holomorphic diffeomorphism $\Phi
:= \tau_{1}\circ\tau_{2}$ of $(\Bbb C^{n+p-1},0)$. 
We assume that it is not tangent to the identity at the origin. 
\begin{defis}
Let ${\cal I}$ be monomial ideal of ${\cal O}_{n+p-1}$.
We shall say that ${\cal I}$ is compatible with $T_1$ and $T_2$  (resp. with the anti-linear involution $\rho$) if the maps $T_j^*:\widehat{\cal O}_{n+p-1}\rightarrow \widehat{\cal O}_{n+p-1}$ (resp. $\rho^*:\widehat{\cal O}_{n+p-1}\rightarrow \overline{\widehat{\cal O}_{n+p-1}}$) defined by $T_j^*(f)=f\circ T_j$ (resp. $\rho^*(f)=f\circ \rho$ ) preserve the splitting $ \widehat{\cal O}_{n+p-1}= \widehat{\cal I}\oplus \widehat{CI}$ (resp. maps  $\widehat{\cal I}$ to $\overline{\widehat{\cal I}}$ and  $\widehat{CI}$ to $\overline{\widehat{CI}}$).
\end{defis}
The main result of this section is the following.
\begin{theos}\label{ideal}
Let $M$ be a third order analytic perturbation of the quadric $\{y_{p+1}=\ldots=y_{n-1}=0, z_n=Q(z',\bar z')\}$, $z'=(z_1,\ldots,z_p)$ as above. Let $\tau_1$, $\tau_2$ be the holomorphic involutions of $(\Bbb C^{n+p-1},0)$ associated to the complexified submanifold ${\cal M}$ of $\Bbb C^{2n}$. Let ${\cal I}$ be a monomial ideal of ${\cal O}_{n+p-1}$ compatible with $D\tau_1(0)$,  $D\tau_2(0)$ and $\rho$.
Assume that the involutions $\tau_1$ and $\tau_2$ are formally linearizable on the ideal ${\cal I}$. Assume furthermore that $D\Phi(0)$ is a diagonal matrix and is diophantine (resp. on ${\cal I}$ whenever ${\cal I}$ is properly embedded). Then the $\tau_j$'s are simultaneously and holomorphically linearizable on ${\cal I}$. Moreover, the linearizing diffeomorphism can be chosen so that it commutes with the linear anti-holomorphic involution $\rho$.
\end{theos}
\begin{proof}
We can apply lemma \ref{rho} to the $(D\Phi(0),D\tau_{1}(0),D\tau_{2}(0))$. We
reorder the vectorspaces so that
$$
\Bbb C^{n+p-1} = E_1\oplus\cdots\oplus E_r\oplus E_{r+1}\oplus\cdots\oplus E_s\oplus G
$$
with $E_1,\ldots, E_h$ are {\it hyperbolic},  $E_{h+1},\ldots, E_r$ are {\it elliptic} and $E_{r+1},\ldots, E_s$ are {\it complex}.
This means that $\mu_i=\bar\mu_i^{-1}$ for $1\leq i\leq h$, $\mu_i=\bar\mu_i$
for $h+1\leq i\leq r$ and $\mu_i\not\in\{\bar\mu_,\bar\mu_i^{-1}\}$ for $r+1\leq i\leq s$. 

Let us consider the ring $\widehat{\cal O}_{n+p-1}^{D\Phi(0)}$ of formal invariants
of $D\Phi(0)$. Since $\mu_i\mu_i^{-1}=1$,  the monomials $\zeta_{i,u}\eta_{i,v}$, $1\leq u,v\leq n(i)$ are
invariants for $1\leq i\leq r$ as well as $\zeta_{i,k,u}\eta_{i,k,v}$,
$k=1,2$, $1\leq u,v\leq n(i)$ and $r+1\leq i\leq s$ and the $\upsilon_j$'s with
$j=1,\ldots, q$ (their associated eigenvalues are $1$). More generally, let $(p,q)\in \Bbb N^{r+2s}\times \Bbb N^{r+2s}$ such that 
$$
\Pi_{i=1}^{r} \mu_{i}^{p_{i}-q_i}\Pi_{j=r+1}^s\mu_{j}^{p_{j}-q_{r}} \bar \mu_{j}^{-p_{j+1}+q_{j+1}}=1.
$$
then all the monomials 
\begin{equation}\label{resonnances}
\Pi_{i=1}^{r}\left(\Pi_{k=1}^{n(i)}\zeta_{i,k}^{p_{i,k}}\eta_{i,k}^{q_{i,k}}\right)\Pi_{j=r+1}^{s}\left(\Pi_{k=1}^{n(j)}\zeta_{j,k}^{p_{j,k}}\eta_{j,k}^{q_{r,k}}\zeta_{j+1,k}^{p_{j+1,k}}\eta_{j+1,k}^{q_{r+1,k}}\right)
\end{equation}
for which $\sum_{k=1}^{n(i)} p_{i,k}=p_i$ and $\sum_{k=1}^{n(i)} q_{i,k}=q_i$ ($1\leq i\leq r+2s$)
are invariants (the $p_{i,k}$'s and $q_{i,k}$'s are non-negative integers). Such a monomial will be written as $\zeta^P\eta^Q$ where
$P=(\{p_{i,k}\}_{1\leq k\leq n(i)})_{1\leq i\leq r+2s}$ and $Q=(\{q_{i,k}\}_{1\leq k\leq n(i)})_{1\leq i\leq r+2s}$. Let ${\cal R}$ denotes the set
of such $(P,Q)$ which generates, together with the $\upsilon_j$'s, the ring of
invariants (in fact, it is a module of finite type; so just the generators can
be selected). 

Let $ResIdeal$ denote the ideal of ${\cal O}_{n-p+1}$ generated by the monomials of the set ${\cal R}$.

Let us write that the $\tau_i$'s are simultaneously formally linearizable on ${\cal I}$.
To be more specific, let us write, for $j=1,2$
$$
\tau_{j}(\zeta',\eta',\upsilon')=\begin{cases}
    \zeta''_i= \lambda_{j,i}\eta'_{i}+ f_{j,i}(\zeta',\eta',\upsilon')\quad
    i=1,\ldots, r+2s\\ \eta''_i=\lambda_{j,i}^{-1}\zeta'_{i}+ g_{j,i}(\zeta',\eta',\upsilon')\quad i=1,\ldots, r+2s\\
\upsilon''= \upsilon'+h_{j}(\zeta',\eta',\upsilon').
\end{cases},
$$
and 
$$
\hat\Psi(\zeta,\eta,\upsilon)=\begin{cases}
\zeta'_i= \zeta_{i}+U_{i}(\zeta,\eta,\upsilon)\quad i=1,\ldots, r+2s\\
\eta'_i=\eta_{i}+V_{i}(\zeta,\eta,\upsilon)\quad i=1,\ldots, r+2s\\
\upsilon'= \upsilon+W(\zeta,\eta,\upsilon).
\end{cases}.
$$
Here we have used to following convention : if $j\geq 1$ then $\mu_{r+2j}=\bar \mu_{r+2j-1}^{-1}$ and $\theta_j:=(\zeta_{r+2j-1},\eta_{r+2j-1},\zeta_{r+2j},\eta_{r+2j})$ are coordinates on $E_{r+j}$.
Let us write that $\hat \Psi$ conjugates $\tau_j$ to $\tilde \tau_j$ with
$$
\tilde \tau_{j}(\zeta,\eta,\upsilon)=\begin{cases}
    \zeta'_i= \lambda_{j,i}\eta_{i}+\tilde f_{j,i}(\zeta,\eta,\upsilon)\quad
    i=1,\ldots, r+2s\\ \eta'_i=\lambda_{j,i}^{-1}\zeta_{i}+ \tilde g_{j,i}(\zeta,\eta,\upsilon)\quad i=1,\ldots, r+2s\\
\upsilon'= \upsilon+\tilde h_{j}(\zeta,\eta,\upsilon).
\end{cases},
$$
where the $\tilde f_{j,i}$'s, the $\tilde g_{j,i}$'s and the components of
$\tilde h_j$ belong to the ideal ${\cal I}$. We have $\hat\Psi\circ \tilde \tau_j= \tau_j \circ \hat
\Psi$; that is
$$
(*)\begin{cases}
\lambda_{j,i}V_i-U_i\circ\tilde \tau_{j}= \tilde f_{j,i}- f_{j,i}\circ\hat\Psi(\zeta,\eta,\upsilon)\quad i=1,\ldots, r+2s\\
\lambda_{j,i}^{-1}U_i-V_i\circ\tilde \tau_{j}=\tilde g_{j,i}-g_{j,i}\circ \hat\Psi(\zeta,\eta,\upsilon)\quad i=1,\ldots, r+2s\\
W-W\circ\tilde \tau_{j}=\tilde h_{j}- h_{j}\circ \hat\Psi.
\end{cases}
$$
Let us find an equation involving only the unkown $U_i$ (resr. $V_i$, $W$). We
recall that $\mu_{i}:= \lambda_{1,i}\lambda_{2,i}^{-1}$. Since we have
$\lambda_{1,i}V_i-U_i\circ \tilde \tau_{1}=\tilde f_{1,i}- f_{1,i}\circ\hat\Psi$, we
have $\lambda_{1,i}V_i\circ \tilde\tau_{2}- U_i\circ \tilde\Phi = \tilde f_{1,i}\circ\tilde
\tau_{2}- f_{1,i}\circ\hat\Psi\circ \tilde \tau_2$. Here, $\tilde \Phi$
denotes $\tilde\tau_1\circ\tilde\tau_2$. According to the second equation of
$(*)$ for $j=2$, we have $V_i\circ \tilde\tau_{2}=\lambda_{2,i}^{-1}U_i-\tilde
g_{2,i}+ g_{2,i}\circ \hat\Psi$. Therefore, we obtain 
$$
\mu_{i}U_i - U_i\circ \tilde\Phi = \left(\tilde f_{1,i}- f_{1,i}\circ\hat\Psi\right)\circ \tilde \tau_2+ \lambda_{1,i}\left(\tilde g_{2,i}- g_{2,i}\circ\hat\Psi\right).
$$
Since the $\tau_j$'s are involutions, we have the following relations~:
%we have 
%\begin{eqnarray*}
%\lambda_{j,i}\left(\lambda_{j,i}^{-1}\zeta'_{i}+ g_{j,i}(\zeta',\eta',\upsilon')\right)+f_{j,i}\circ \tau_j & = & \zeta'_{i}\\
%\lambda_{j,i}^{-1}\left(\lambda_{j,i}\eta'_{i}+ f_{j,i}(\zeta',\eta',\upsilon')\right)+g_{j,i}\circ \tau_j & = & \eta'_{i}\\
%h_j\circ\tau_j + h_j &= & 0;
%\end{eqnarray*}
%that is 
\begin{eqnarray}
\lambda_{j,i}g_{j,i}+f_{j,i}\circ \tau_j & = & 0\label{rel1}\\
\lambda_{j,i}^{-1}f_{j,i}+g_{j,i}\circ \tau_j & = & 0\label{rel3}\\
h_j\circ\tau_j + h_j &= & 0\label{rel2}.
\end{eqnarray}
We have the following relations among the $f$'s and the $g$'s~:
\begin{equation}\label{f-g}
g_{1,i}\circ \tau_{2} - \mu_{i}^{-1} g_{2,i}\circ \tau_{2}= -\lambda_{1,i}^{-1}\left(f_{1,i}\circ \tau_{1} - f_{2,i}\circ \tau_{2}\right)\circ \tau_{2}.
\end{equation}
Using the fact that $f_{1,i}\circ\hat\Psi\circ \tilde \tau_{2}= f_{1,i}\circ \tau_2\circ \hat\Psi$, we find the following relations~:
\begin{equation}\label{u}
\mu_{i}U_i - U_i\circ \tilde\Phi =  \left(\tilde f_{1,i}- \mu_{i}\tilde f_{2,i}\right)\circ \tilde\tau_{2} - \left( f_{1,i}- \mu_{i} f_{2,i}\right)\circ \tau_{2}\circ\hat\Psi=:\alpha_i.
\end{equation}
Similarly, we obtain 
\begin{equation}\label{v}
\mu_{i}^{-1}V_i - V_i\circ \tilde\Phi =  \left(\tilde g_{1,i}- \mu_{i}^{-1}\tilde g_{2,i}\right)\circ \tilde\tau_{2} - \left( g_{1,i}- \mu_{i}^{-1} g_{2,i}\right)\circ \tau_{2}\circ\hat\Psi=:\beta_i,
\end{equation}
as well as 
\begin{equation}\label{w}
W- W\circ \tilde\Phi =  \left(\tilde h_{1}\circ \tilde\tau_{2}+\tilde h_{2}\right)- \left(h_{1}\circ\tau_2+h_2\right)\circ\hat\Psi=:\gamma.
\end{equation}

The germ of diffeomorphism $\Phi$ is formally linearizable on ${\cal I}$. If
its linear part $D\Phi(0)$ is diophantine (resp. on  ${\cal I}$ if ${\cal I}$ is properly embedded), then, according to
theorem \ref{theo-invariant}, $\Phi$ is actually holomorphically linearizable
on ${\cal I}$ and by a unique normalizing transformation which projection
on the vector space spanned by $\widehat{\cal I}$ and the formal centralizer
of $D\Phi(0)x$ is zero. Let $\Psi'$ be precisely this germ of
holomorphic diffeomorphism of $(\Bbb C^{n+p-1},0)$~:
$$
\Psi'(\zeta,\eta,\upsilon)=\begin{cases}
\zeta'_i= \zeta_{i}+U'_{i}(\zeta,\eta,\upsilon)\quad i=1,\ldots, r+2s\\
\eta'_i=\eta_{i}+V'_{i}(\zeta,\eta,\upsilon)\quad i=1,\ldots, r+2s\\
\upsilon'= \upsilon+W'(\zeta,\eta,\upsilon).
\end{cases}
$$
Let us show that it commutes with $\rho$. First of all, we have $\rho\Phi\rho=\Phi^{-1}$. Then, according to the second part of theorem \ref{theo-invariant}, it is sufficient to show that $\rho{\cal C}_{D\Phi(0)x}\rho = {\cal C}_{D\Phi(0)x}$. Let $R\in {\cal C}_{D\Phi(0)x}$. We have $D\Phi(0)R(x) = R D\Phi(0)(x)$. We recall that $D\Phi(0)x=T_1T_2x$. Therefore, we have
$$
T_1T_2\rho R\rho = T_1\rho T_1 R\rho = \rho T_2 T_1 R\rho = \rho R T_2 T_1 \rho=\rho R\rho T_1T_2.
$$

Let us modify $\Psi'$ in such a way that it linearizes simultaneously
$\tau_1,\tau_2$ on ${\cal I}$. Let $pr(f)$ (resp. $pr_{\cal I}(f)$) denote the projection of the formal
power series $f$ on the vector space $\widehat{CI}$ spanned by the $x^Q$'s which doesn't belong
to $\widehat{\cal I}$ (resp. on $\widehat{\cal I}$). Let us project the equations $(*)$ onto this space. We find
$$
(**)\begin{cases} \lambda_{j,i}pr(V_i)-pr(U_i\circ T_{j})= - pr(f_{j,i}\circ\hat\Psi)\quad i=1,\ldots, r+2s\\
\lambda_{j,i}^{-1}pr(U_i)-pr(V_i\circ T_{j})=-pr(g_{j,i}\circ \hat\Psi)\quad i=1,\ldots, r+2s\\
pr(W)-pr(W\circ T_{j})=- pr(h_{j}\circ \hat\Psi).
\end{cases}
$$
On the other hand, projecting the equations $(\ref{u})$,$(\ref{v})$ and $(\ref{w})$ on $CI$ we obtain
\begin{eqnarray*}
\mu_{i}pr(U_i) - pr(U_i)\circ D & = &  - pr\left(\left( f_{1,i}- \mu_{i} f_{2,i}\right)\circ \tau_{2}\circ\hat\Psi\right)\\
\mu_{i}^{-1}pr(V_i) - pr(V_i)\circ D & = &  - pr\left(\left( g_{1,i}- \mu_{i}^{-1} g_{2,i}\right)\circ \tau_{2}\circ\hat\Psi\right)\\
pr(W)- pr(W\circ D) &= & - pr\left(\left(h_{1}\circ\tau_2+h_2\right)\circ\hat\Psi\right).
\end{eqnarray*}

Let ${\cal V}_{\pm i}$ (resp. ${\cal V}_0$), $1\leq i\leq r$, be the closed
subspace of $\Bbb C[[\zeta,\eta,\upsilon]]$ generated by monomials
$\zeta^{Q_1}\eta^{Q_2}\upsilon^{Q_3}$ for which
$\mu^{Q_1-Q_2}=\mu_{i}^{\pm 1}$ (resp. $\mu^{Q_1-Q_2}=1$). 
Let $P_i$ be the associated projection. Applying these various projections to the previous equations leads to the following compatibility equations~:
\begin{eqnarray}
P_i\left(pr\left(\left( f_{1,i}- \mu_{i} f_{2,i}\right)\circ \tau_{2}\circ\hat\Psi\right)\right)&=&0\nonumber\\
P_{-i}\left(pr\left(\left( g_{1,i}- \mu_{i}^{-1} g_{2,i}\right)\circ \tau_{2}\circ\hat\Psi\right)\right)&=&0\label{compat-g}\\
P_0\left(pr\left(\left(h_{1}\circ\tau_2+h_2\right)\circ\hat\Psi\right)\right) &= & 0\nonumber.
\end{eqnarray}

According to the properties of the eigenvalues of $D\Phi(0)$ and its invariants, we have the following~:
\begin{itemize}
\item If $Q:=(Q_1,Q_2,Q_3)\in\Bbb N^{p+p+q}$ is such that $\mu^{Q_1-Q_2}=1$ then $\zeta^{Q_1}\eta^{Q_2}\upsilon^{Q_3}\circ T_{j}$ is a monomial $\zeta^{Q'_1}\eta^{Q'_2}\upsilon^{Q'_3}$ such that $\mu^{Q'_1-Q'_2}=1$.\\
\item If $Q:=(Q_1,Q_2,Q_3)\in \Bbb N^{2p+q}$ and $1\leq i\leq p$ is such that
  $\mu^{Q_1-Q_2}=\mu_{i}^{\pm 1}$ then
  $\zeta^{Q_1}\eta^{Q_2}\upsilon^{Q_3}\circ T_{j}$ is a monomial
  $\zeta^{Q'_1}\eta^{Q'_2}\upsilon^{Q'_3}$ such that
  $\mu_k^{Q'_1-Q'_2}=\mu_{i}^{\mp 1}$.\\
\end{itemize}

Since ${\cal I}$ is compatible with the $T_j$'s, we have 
$pr(f\circ T_j)=pr(f)\circ T_j$ for any formal power series $f$. Moreover, we
have
$$
P_{\pm i}(f\circ T_j)= P_{\mp i}(f)\circ T_j,\quad P_{0}(f\circ T_j)= P_{0}(f).
$$
We have some freedom in the choice of the normalizing transformation of the
$\tau_j$'s. Let us show that there is a unique solution $U_i$, $V_i$ and $W$ such that 
$$
pr_{\cal I} U_i=pr_{\cal I} V_i=pr_{\cal I} W=0,
$$
and 
\begin{equation}\label{pnormal}
P_{-i}(pr(V_i))=0;\quad P_0(pr(W))=0.
\end{equation}
\begin{rems}
If the formal centralizer $\widehat{\cal C}_{D\phi(0)x}$ is contained in $\widehat{\cal I}$ then the last condition is always satisfied.
\end{rems}
In fact, let apply $P_{-i}$ (resp. $P_{i}$, $P_0$) to the first (resp. second,
last)
equation of $(**)$ and let us show that system obtained is solvable. We have
\begin{eqnarray*}
P_i(pr(U_i))\circ T_j & = & P_{-i}(pr(f_{j,i}\circ\hat\Psi))\\
\lambda_{j,i}^{-1}P_i(pr(U_i)) & = & -P_{i}(pr(g_{j,i}\circ\hat\Psi))\\
P_0(pr(h_j\circ\hat \Psi)) & = & 0.
\end{eqnarray*}

The last equation is obtained from equation $(\ref{rel2})$. In fact, we have
$$
0=h_j\circ\tau_j\circ\hat\Psi+h_j\circ\hat\Psi =(h_j\circ\hat\Psi)\circ\tilde\tau_j+(h_j\circ\hat\Psi).
$$
Since $\tilde \tau_j$ is linear on ${\cal I}$, we have 
$$
pr((h_j\circ\hat\Psi)\circ\tilde\tau_j)=pr((h_j\circ\hat\Psi)\circ
T_j)=pr(h_j\circ\hat\Psi)\circ T_j.
$$
Therefore, we have $P_0(pr(h_j\circ\hat\Psi))=0$.
The compatibility condition of the the first equation is 
$$
P_{-i}(pr(f_{j,i}\circ\hat\Psi))\circ T_j+\lambda_{j,i}P_{i}(pr(g_{j,i}\circ\hat\Psi))=0.
$$
Let us show that it can be obtained from equation $(\ref{rel1})$. In fact, let us
first compose on the right by $\hat\Psi$ and then project onto $CI$. We have
$$
\lambda_{j,i}pr(g_{j,i}\circ\hat\Psi)+pr(f_{i,j}\circ\hat\Psi\circ \tilde\tau_j ) =  0.
$$
As above, this leads to 
$$
\lambda_{j,i}pr(g_{j,i}\circ\hat\Psi)+pr(f_{i,j}\circ\hat\Psi\circ T_j ) =  0,
$$
and then to
$$
\lambda_{j,i}P_i\left(pr(g_{j,i}\circ\hat\Psi)\right)+P_{-i}\left(pr(f_{i,j}\circ\hat\Psi)\right)\circ T_j  =  0.
$$

Let us consider the following equations
\begin{equation}
u_i  =  -\lambda_{j,i}P_{i}(pr(g_{j,i}(\zeta+U'+u,\eta+V',\upsilon +
W'))),\quad i=1,\ldots, r+2s,\label{equ-u}
\end{equation}
where the $u_i$'s are the unknowns and which components belong to the range of $P_i\circ pr$. Here, the maps $U'$, $V'$ and $W'$ are
holomorphic maps defined by the $\Psi'$, the normalized transformation which linearizes $\Phi$ on
${\cal I}$. Since the $g_{j,i}$'s are non-linear holomorphic functions, we can
apply the implicit function theorem to obtain holomorphic $u_i$'s in a
neighbourhood of the origin satisfying equation $(\ref{equ-u})$. The map $u$ is alspo solution of the equations $u_i\circ T_1  =  P_{-i}(pr(f_{1,i}\circ (\Psi'+u)))$. Then, the
germ of holomorphic diffeomorphism at the origin of $\Bbb C^{r+n-1}$ defined
to be
$$
\Psi(\zeta,\eta,\upsilon)=
\begin{cases}
\zeta'_i= \zeta_{i}+U'_{i}(\zeta,\eta,\upsilon)+u_i(\zeta,\eta,\upsilon),\quad i=1,\ldots, r+2s\\
\eta'_i=\eta_{i}+V'_{i}(\zeta,\eta,\upsilon),\quad i=1,\ldots, r+2s\\
\upsilon'= \upsilon+W'(\zeta,\eta,\upsilon)
\end{cases}
$$
linearizes simultaneously and holomorphically the $\tau_j$'s on ${\cal I}$. It is the unique germ of diffeomorphism which linearizes the $\tau_j$'s and satisfies to $(\ref{pnormal})$ and $(\ref{equ-u})$. It remains to show that$\Psi$ commutes whith the anti-holomorphic involution $\rho$. In fact, we have $\hat\Psi\circ \tilde\tau_2 = \tau_2\circ \hat\Psi$. By composition on the right and on the left by $\rho$, we obtain
$$
\left(\Psi'+(\rho\circ u\circ\rho)\right)\circ (\rho\circ\tilde \tau_2\circ\rho) = \tau_{1}\circ \left(\Psi'+(\rho\circ u\circ\rho)\right).
$$
This is due to the fact that $\rho$ is an involution which commutes with $\Psi'$. It is to be noticed that the map $\rho\circ\tilde \tau_2\circ\rho$ is equal to $T_1$ modulo the ideal ${\cal I}$. Therfore, if we project outside ${\cal I}$, we obtain
$$
pr(\Psi')\circ T_1+pr (\rho\circ u\circ\rho)\circ T_1 = pr\left(\tau_{1}\circ \left(\Psi'+(\rho\circ u\circ\rho)\right)\right).
$$
According to the properties of $\Psi'$, we have obtain
$$
P_i(pr (\rho\circ u\circ\rho)_i)\circ T_1 = P_{-i}(pr\left(f_{1,i}\circ \left(\Psi'+(\rho\circ u\circ\rho)\right)\right)).
$$
Therefore, $\rho\circ u\circ\rho$ is solution of the same equation as $u$. By uniqueness, we obtain $\rho\circ u\circ\rho=u$. Hence, $\Psi$ commutes with $\rho$. This end the proof of the theorem.
\end{proof}

%Let us consider the first equation of $(**)$. Let us $u$ the vector $(P_i(pr(U_i)))_{i=1+r+2s}$. It is solution of the equation
%$$
%u T_j= (P_{-i}(pr(f_{j,i}(\Psi'+u))))_{i=1+r+2s}.
%$$
%By composition on the right and on the left by $\rho$, we obtain
%$$
%\rho u\rho T_{j-1}= \rho(P_{-i}(pr( f_{j,i}\rho(\rho\Psi'+\rho u))))_{i=1+r+2s}\rho.
%$$

Let us express a direct geometric implication. 
%First of all, we have the
\begin{coros}
Under the assumptions of the theorem and in the new holomorphic coordinates system, the complexified submanifold ${\cal M}$ intersects each irreducible component of the zero locus $V(\cal I)$ which is invariant under $\rho$, $T_1$ and $T_2$   along a complex submanifold which is a complexified quadric of $\Bbb C^k$, $k\leq n$ (if the $T_j$'s are still involutions) .
\end{coros}
\begin{proof}
Under the assumptions of the theorem, the $\tau_j$'s are holomorphically linearizable on ${\cal I}$. Hence, in the new coordinates, their restriction to any irreducible component of $V(\cal I)$ is equal to their restriction of their linear part at the origin. Assume that $V_i$ is an irreductible component of $V({\cal I})$ which is invariant under $\rho$, $T_1$ and $T_2$.
Since ${\cal I}$ is a monomial ideal, its zero locus is an intersection of unions of complex hyperplanes $\{\zeta_{i_j}=0\}$, $\{\eta_{i_j}=0\}$ and $\{\upsilon_{k_l}=0\}$. Hence it is a linear complex manifold.
Then, the submanifold $V_i\cap Fix(\rho)$ of $M$ is a quadric associated to the linear involutions $T_{j|V_j}$.
%. 
%Let $V_i$ be an irreductible component of $V({\cal I})$ invariant under $\rho$. The restrictions of the $\tau_j$'s to $V_i$ (which are linear) are the involutions of a quadric of dimension less or equal to $n$.
\end{proof}

\subsection{Holomorphic equivalence to quadrics}

As a direct consequence, we have the following result~:
\begin{theos}
Assume that $M$ is formally equivalent to the quadric defined by its quadratic
part. Let $\Phi$ be the associated germ of holomorphic diffeomorphism. Assume
that its linear part $D\Phi(0)$ is semi-simple and diophantine, then $M$ is biholomophic to the quadric
\begin{equation} (Q)
\begin{cases} z_1 = x_1+iy_1\\ \vdots \\ z_p = x_p+iy_p\\ 
y_{p+1} =  0\\ \vdots\\ y_{n-1} =  0\\
z_{n} = Q(z',\bar z')\end{cases}
\end{equation}
\end{theos}
When $p=1$, $n=2$, this result is due to X. Gong
\cite{gong-hyperbolic} in the hyperbolic case.
\begin{proof}
Let us apply the main theorem with the zero ideal ${\cal I}:=(0)$. In fact, it is compatible with the $D\tau_j(0)$'s and $\rho$. That $M$ is formally equivalent to the quadric $(Q)$ means precisely that the $\tau_j$'s are simultaneously formally linearizable. According to the assumption, the $\tau_j$'s are simultaneously  and holomorphically linearizable. Since the linearizing diffeomorphism $\Psi$ commutes with the linear anti-holomorphic involution $\rho$, $\Psi$ is the complexified of an holomorphic diffeomorphism $\psi$ of $(\Bbb C^n,0)$ which maps $M$ to the quadric $(Q)$.
\end{proof}

\subsection{Cutting varieties}

We want to apply the theorem to the ideal ${\cal I}:= ResIdeal$, the ideal generated by the set ${\cal R}$ of resonnant monomials $(\ref{resonnances})$. Almost all pair of involution $(\tau_1,\tau_2)$ are formally linearizable along the resonant ideal ${\cal I}$. In fact, the normal form theory for the diffeomorphism $\Phi$ tells us that all but a finite number of the monomials appearing in the Taylor expansion of a normal form belongs to the resonannt ideal. 
\begin{lemms}
Assume that the non-linear centralizer ${\cal C}_{D\Phi(0)}$ is included in the resonnant ideal $ResIdeal$ and that $D\Phi(0)$ has distinct eigenvalues. Then the $\tau_j$'s are formally and simultaneously linearizable on ${\cal I}=ResIdeal$. Moreover, ${\cal I}$ is compatible with the $D\tau_j(0)$'s and $\rho$.
\end{lemms}
\begin{proof}
Let us write 
$$
\Phi=
\begin{cases}
\mu_i\zeta_i+\phi_i, \quad i=1,\ldots, r+2s\\
\mu_i^{-1}\eta_i+\psi_i,\quad i=1,\ldots, r+2s\\
\upsilon+\theta.
\end{cases}
$$
Since $\Phi=\tau_1 \circ \tau_2$, we have
$$
\begin{cases}
\lambda_{1,i}g_{2,i}+f_{1,i}\circ\tau_2= \phi_i\\
\lambda_{1,i}^{-1}f_{2,i}+g_{1,i}\circ\tau_2= \psi_i\\
h_{2}+h_{1}\circ\tau_2= \theta.
\end{cases}
$$
But according to relations $(\ref{rel1})$, $(\ref{rel3})$ and $(\ref{rel2})$, we have
\begin{equation}\label{rel4}
\begin{cases}
g_{2,i}= -\lambda_{2,i}^{-1}f_{2,i}\circ\tau_2 \\
g_{1,i}= -\lambda_{1,i}^{-1}f_{1,i}\circ\tau_2 \\
h_2+h_2\circ \tau_2=0.
\end{cases}
\end{equation}
Therefore, we have
$$
\begin{cases}
-(\lambda_{1,i}\lambda_{2,i}^{-1})f_{2,i}\circ\tau_2+f_{1,i}\circ\tau_2= \phi_i\\
\lambda_{1,i}^{-1}f_{2,i}-\lambda_{1,i}^{-1}f_{1,i}\circ \Phi= \psi_i\\
-h_{2}\circ \tau_1+h_{1}\circ\tau_2= \theta.
\end{cases}
$$
We find that
\begin{equation}\label{rel5}
\begin{cases}
\mu_i^{-1}f_{1,i}-f_{1,i}\circ \Phi = \mu_{i}\psi_i+\mu_{i}^{-1}\phi_i\circ\tau_2\\
f_{2,i}=f_{1,i}\circ \Phi+ \mu_i\psi_i
\end{cases}
\end{equation}
Assume that the $\phi_i$'s, $\psi_i$'s and $\theta$ belong to ${\cal I}$. Assume that the $k$-jet of $\tau_1-T_1$ and $\tau_2-T_2$ belong to also the ideal ${\cal I}$. Then, by equalities $(\ref{rel4})$ and $(\ref{rel5})$, we have that $\mu_i^{-1}f_{1,i,k+1}-f_{1,i,k+1}\circ D\Phi(0)$ belong to ${\cal I}$. Here, $f_{1,i,k+1}$ denotes the homogenous polynomial of order $k+1$ in the Taylor expansion of $f_{1,i}$ at the origin. Since the centralizer of $D\Phi(0)$ is contained in ${\cal I}$ then $f_{1,i,k+1}$ belongs to ${\cal I}$. It follows that $f_{2,i, k+1}$ also belongs to ${\cal I}$. So, are the $g_{j,i,k+1}$'s and the $h_{j,k+1}$'s. Therefore, under the assuptions, the $\tau_i$'s are formally linearizable on the ideal ${\cal I}$ since $\Phi$ is.

About the second point, the $D\tau_j(0)$'s leave each monomial $\zeta_i\eta_i$ invariant. Moreover, if a monomial $\zeta^Q\eta^R$ is first integral of $D\Phi(0)$, then so is $\eta^Q\zeta^R$. this is due to the fact that if $\mu^Q(\mu^{-1})^R=1$ then $(\mu^{-1})^Q(\mu)^R=1$. As a consequence, each first integral of $D\Phi(0)$ is sent to a first integral of $D\Phi(0)$ by the $D\tau_j(0)$'s. Furthermore, $\bar\rho^*$ maps a monomial first integral to another. This is due to the fact that the eigenvalues are distinct.
\end{proof}
The assumption of the lemma means that the resonnances are generated by the generators of the ring of invariants of the linear part of $\Phi$ at the origin.
\begin{theos}
Assume that the non-linear centralizer ${\cal C}_{D\Phi(0)x}$ is included in the resonnant ideal $ResIdeal$. Moreover let us assume that $D\Phi(0)$ has distinct hyperbolic eigenvalues and is diophantine (on $ResIdeal$ if the latter is properly embedded). Then, there are good holomorphic coordinates $(v_1,\ldots,v_n)$ of $\Bbb C^n$, in which the submanifold $M$ intersects the complex linear manifold $\{v_{p+1}=\ldots =v_n=0\}$ along a real analytic subset $V$ whose complexified is nothing but the zero locus of the resonnant ideal $ResIdeal$~: $V=V(ResIdeal)\cap Fix(\rho)$, where $Fix(\rho)$ denotes the fixed point set of $\rho$.
\end{theos}
\begin{proof}
%Let us show that both $Q(z',w')$ and $\bar Q(w',z')$ vanish on the zero locus $V(ResIdeal)$ of $ResIdeal$. 
First of all, we make a change of coordinates as in lemma \ref{rho}. Then, the quadratic functions which are invariant under both $T_1$ and $T_2$ are linear combinations of monomials of the form $\zeta_i\eta_i$ and $\upsilon_j\upsilon_k$. Therefore, in these coordinates, the functions defining the complexified quadric are functions $f,g$ of the $\zeta_i\eta_i$'s and $\upsilon_j\upsilon_k$'s. Since these monomials vanish on the zero locus of $ResIdeal$, so does $f$ and $g$. We apply theorem \ref{ideal}~: $V(ResIdeal)$ is invariant under the $\tau_j$'s and their restriction to it are equal to $T_{j|V(ResIdeal)}$. Because of hyperbolicity, each irreductible component of $V(ResIdeal)$ is invariant under $\rho$. Hence, the functions which are invariants by the $\tau_j$'s vanish on zero on the zero locus of $Resideal$.
\end{proof}
When $p=1$, $n\geq 2$ and in the hyperbolic case, this result is due to
Wilhelm Klingenberg Jr. \cite{klingenberg}. He found, for instance in the case $n=2$, that, in some suitable holomorphic coordinates, the complex analytic set $\{v_2=0\}$ intersects $M$ along
$\{(\zeta_1,\eta_1)\in \Bbb R^2,\;\zeta_1\eta_1=0\}$. In fact, in this case, the ring of formal invariants of $\Phi$ is generated by only one element: $\zeta_1\eta_1$. 

\begin{rems}
When the spectrum of $D\Phi(0)$ contains {\bf elliptic} eigenvalues and {\bf complex} eigenvalues, then one has to use the set 
$$
\tilde V=V({\cal I})\cap \left(\bigcap_{i\in {\cal E}} \{\zeta_i=\eta_i=0\}\right))\cap \left(\bigcap_{j\in {\cal C}} \{\zeta_{r+2j-1}=\eta_{r+2j-1}=\zeta_{r+2j}=\eta_{r+2j}=0\}\right)
$$
which is invariant by $\rho$. Then, the conclusion of the previous theorem holds with $\tilde V\cap Fix(\rho)$.
\end{rems}

%%% Local Variables: 
%%% mode: latex
%%% TeX-master: "webster"
%%% TeX-master: "webster"
%%% End: 

\bibliographystyle{alpha}
\bibliography{normal,math,asympt,analyse,stolo,cr}
\end{document}